\documentclass[10pt,onecolumn,twoside,draftclsnofoot]{IEEEtran}
\usepackage{url}

\usepackage{booktabs} 
\usepackage{array} 
\usepackage{paralist} 
\usepackage{verbatim} 
\usepackage{alltt}
\usepackage{amsmath}
\usepackage[hyperfootnotes=false]{hyperref}
\usepackage{subcaption}
\usepackage{graphicx} 
\usepackage{dsfont}
\usepackage{mathtools}
\usepackage{cleveref}
\usepackage{cite}
\usepackage{amssymb}
\usepackage{color}

\newcommand{\qed}{\hfill \ensuremath{\blacksquare}}
\newcommand{\changed}[1]{#1}

\newtheorem{theorem}{Theorem}

\newtheorem{lemma}{Lemma}

\newtheorem{remark}{Remark}

\newcommand{\1}{\mathbf{1}}
\newcommand{\0}{\mathbf{0}}
\newcommand{\kp}{\kappa}

\newcommand{\enma}[1]   {\ensuremath{#1}}

\newcommand{\diag}  {\enma{\mathrm{diag}}}

\newcommand{\trace} {\enma{\mathrm{trace}}}

\newcommand{\blockdiag} {\enma{\mathrm{blockdiag}}}

\title{Skewless Network Clock Synchronization Without Discontinuity: Convergence and Performance}

\author{\IEEEauthorblockN{Enrique Mallada$^*$, Xiaoqiao Meng$^\dagger$, Michel Hack$^\dagger$, Li Zhang$^\dagger$, and Ao Tang$^\flat$}
\\[1ex]
\IEEEauthorblockA{$^*$ Computational and Mathematical Sciences, Caltech, Pasadena, CA 91125, USA.
\\
$^\dagger$ IBM T. J. Watson Research Center. 1101 Kitchawan Rd, Yorktown Heights, NY 10598, USA.\\
$^\flat$ School of ECE, Cornell University, Ithaca, NY 14853, USA.
}
\thanks{Preliminary version of this paper appears in IEEE ICNP 2013 \cite{Mallada:2013hg}.}
}

\begin{document}

\maketitle

\begin{abstract}
This paper examines synchronization of computer clocks connected via a data network and proposes a skewless algorithm to synchronize them. Unlike existing solutions, which either estimate and compensate the frequency difference (skew) among clocks or introduce offset corrections that can generate jitter and possibly even backward jumps, our solution achieves synchronization without these problems. We first analyze the convergence property of the algorithm and provide explicit necessary and sufficient conditions on the parameters to guarantee synchronization. We then study the effect of noisy measurements (jitter) and frequency drift (wander) on the offsets and synchronization frequency, and further optimize the parameter values to minimize their variance. Our study reveals a few insights, for example, we show that our algorithm can converge even in the presence of timing loops and noise, provided that there is a well defined leader. This marks a clear contrast with current standards such as NTP and PTP, where timing loops are specifically avoided. Furthermore, timing loops can even be beneficial in our scheme as it is demonstrated that highly connected subnetworks can collectively outperform individual clients when the time source has large jitter. The results are supported by experiments running on a cluster of IBM BladeCenter servers with Linux.
\end{abstract}

%

\section{Introduction}\label{sec:introduction}

Keeping consistent time among different nodes in a network is a fundamental requirement of many distributed applications. Nodes' internal clocks are usually not accurate enough and tend to drift  apart from each other over time, generating inconsistent time values. Network clock synchronization allows these devices to correct their clocks to match a global reference of time, such as the Universal Coordinated Time (UTC), by performing  time measurements through a network. For example, for the Internet, network clock synchronization has been an important subject of research and several different protocols have been proposed~\cite{mills_network_2010, sobeih_almost_2007, Anonymous:2008dl, Veitch:2009td, Carli:2014gd, Mallada:2011iu, Ridoux:2012ch}.
These protocols are used in various applications with diverse precision requirements such as banking, communications, traffic monitoring and security. In modern wireless cellular networks, for instance, time-sharing protocols need an accuracy of several microseconds to guarantee the efficient use of channel capacity. Another example is the recently announced Google Spanner~\cite{corbett_spanner_2012}, a globally-distributed database, which depends on globally-synchronized clocks within at most several milliseconds drifts.

The current {\it de facto} standard for IP networks is the Network Time Protocol (NTP) proposed by David Mills~\cite{mills_network_2010}. It is a low-cost, purely software-based
solution whose accuracy mostly ranges from hundreds of microseconds to several milliseconds. On the other hand,
IEEE 1588 (PTP)~\cite{Anonymous:2008dl}  gives superior performance by achieving  sub-microsecond or even nanosecond accuracy. However, it is relatively expensive as it requires special hardware support to achieve those accuracy levels and may not be fully compatible with legacy cluster systems.

Newer synchronization protocols have been proposed with the objective of balancing between accuracy and cost. For example,
 IBM Coordinated Cluster Time (CCT)~\cite{froehlich_achieving_2008} is able to provide better performance than NTP without additional hardware. Its success is based on a skew estimation mechanism~\cite{zhang_clock_2002} that progressively adapts the clock frequency without offset corrections.
\changed{Another solution that  achieves this objective is the RADclock~\cite{Veitch:2009td,Ridoux:2012ch} which decouples skew compensation from offset corrections by decomposing the clock into a high performance difference clock for measuring time differences and a less precise absolute clock that provides UTC time.}



There are two major difficulties that make the network clock synchronization problem challenging.  Firstly, the frequency of hardware clocks is sensitive to temperature and is constantly varying. Secondly, the latency introduced by the OS and network congestion delay results in errors in the time measurements which can be propagated through the network. Thus, most protocols introduce different ways of estimating the frequency mismatch (skew)\cite{zhang_clock_2002, Kim:2012db} and measuring the time difference (offset) \cite{Elson:2002ga,hunt_network_2010} while maintaining a simple network topology~\cite{mills_network_2010, Anonymous:2008dl}.
This leads in particular to extensive literature on skew estimation \cite{Kim:2012db,Marouani:2008fe,moon_estimation_1999,lemmon_model-based_2000}  which suggests that explicit skew estimation is necessary for clock synchronization.

\changed{This paper takes a different approach and shows that using skew estimation is unnecessary.} We provide a simple algorithm that is able to compensate the clock skew without any explicit estimation of it. Our algorithm only uses current offset information and an exponential average of the past offsets. Thus, it neither needs to store long offset history nor perform expensive computations on them. \changed{The solution provided in this paper achieves microsecond level accuracy without requiring any special hardware. Since we do not explicitly estimate the skew, the implementation is simpler and more robust to noise than IBM CCT, and does not introduce offset corrections, which avoids the need of decomposing the clock into several components to reduce jitter as in RADclock.}

By looking at the synchronization problem from a new angle, this paper also provides several new insights. For example, a common practice in the clock synchronization community is to avoid timing loops in the network~ \cite[p. 6]{mills_network_2010} \cite[p. 16, s. 6.2]{Anonymous:2008dl}.  This is because it is thought that timing loops can introduce instability as stated in \cite{mills_network_2010}: {\it
"Drawing from the experience of the telephone industry, which learned such lessons at considerable cost, the subnet topology...
 must never be allowed to form a loop."
} Even though for some parameter values loops can produce instability, we show that a set of proper parameters can guarantee convergence even in the presence of loops. Furthermore, we experimentally demonstrate in Section \ref{sec:experimental_results} that timing loops among clients can actually help reduce the jitter of the synchronization error and is therefore desirable.

The rest of the paper is organized as follows. In Section \ref{sec:comp_clocks} we provide some background on how clocks are actually implemented in computers and how different protocols discipline them. Section \ref{sec:algorithm} motivates and describes our algorithm together with an intuitive explanation of why it works. In Section \ref{sec:analysis}, we analyze the convergence property of the algorithm and determine the set of parameter values and connectivity patterns under which synchronization is guaranteed. The parameter values that guarantee synchronization depend on the network topology, but there exists a subset of them that is independent of topology and therefore of great practical interest. The effect of noisy measurement and wander is studied in Section \ref{sec:noise}, together with an optimization procedure that finds optimal parameter values. Experimental results evaluating the performance of the algorithm are presented in Section \ref{sec:experimental_results}. We conclude in Section \ref{sec:conclusions}.

\section{Computer Clocks and Synchronization} \label{sec:comp_clocks}

Most computer architectures keep their own estimate of time using a counter that is periodically increased by either hardware or kernel's interrupt service routines (ISRs). On Linux platforms for instance, there are usually several different clock devices that can be selected as the clock source by changing the $clocksource$ kernel parameter.
One particular counter that has recently been used by several clock synchronization protocols~\cite{froehlich_achieving_2008,Veitch:2009td} is the Time Stamp Counter (TSC) that counts the number of CPU cycles since the last restart of the system. For example, in the IBM BladeCenter LS21 servers, the TSC is a $64$-bit counter that increments every $\delta^o = 0.416$ns since the CPU nominal frequency $f^o=1/\delta^o=2399.711$MHz.

Based on this counter, each server builds its own estimate $x_i(t)$ of the global time reference, UTC, denoted here by $t$. For example, if $c_i(t)$ denotes the counter's value of computer $i$ at time $t$, then $x_i(t)$ can be computed using
\begin{equation}\label{eq:TSC_map}
x_i(t) = \delta^o c_i(t) +x_i^o,
\end{equation}
where $x^o$ is the estimate of the time when the server was turned on ($t_0$).

\begin{figure}[htp]
       \begin{subfigure}[b]{0.49\columnwidth}
               \centering
               \includegraphics[width=.95\columnwidth,height=.65\columnwidth]{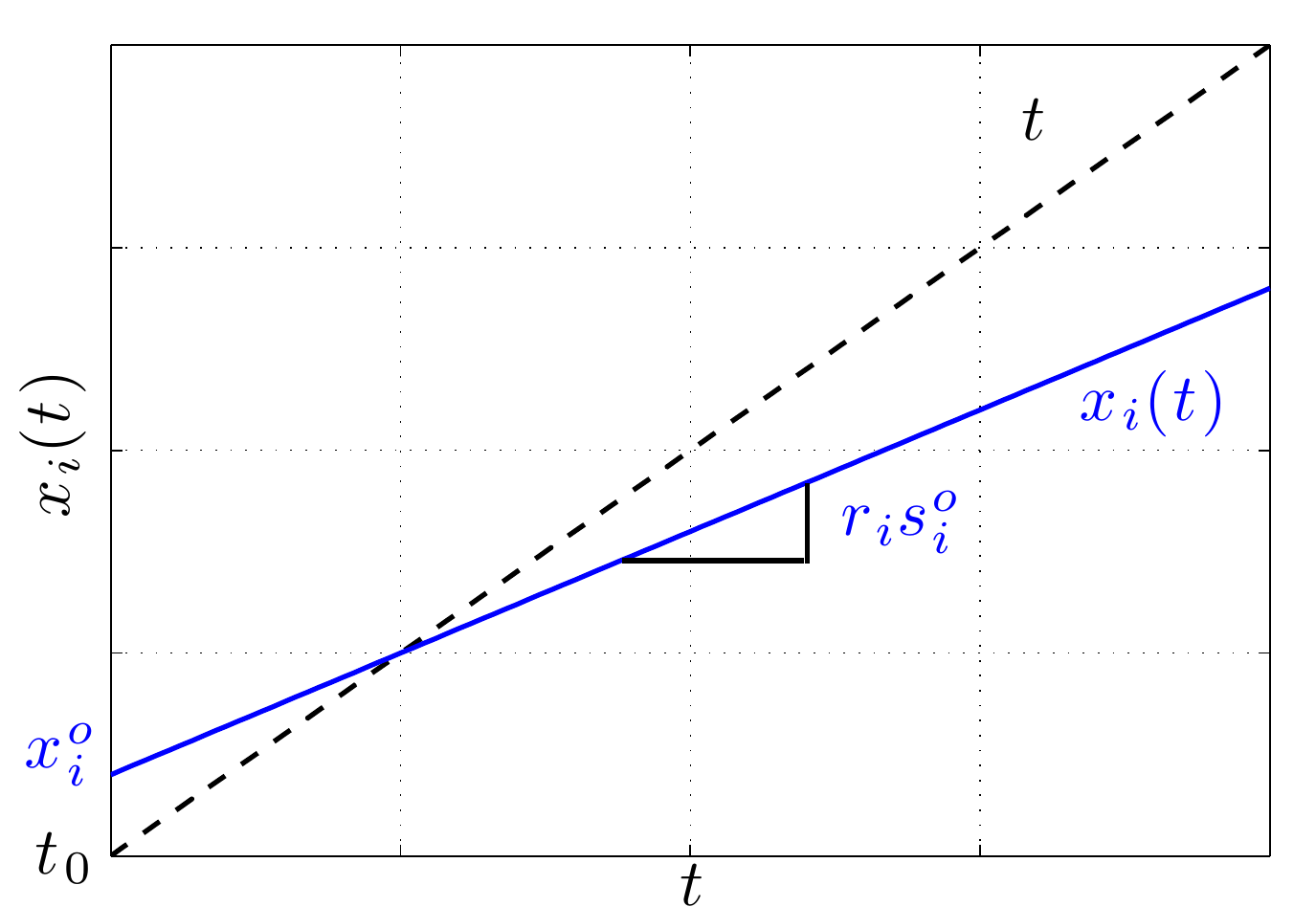}
               \caption{Illustration of computer time estimate $x_i(t)$ and UTC time $t$}\label{fig:xt}
        \end{subfigure}
        \begin{subfigure}[b]{0.49\columnwidth}
               \centering
               \includegraphics[width=.95\columnwidth,height=.65\columnwidth]{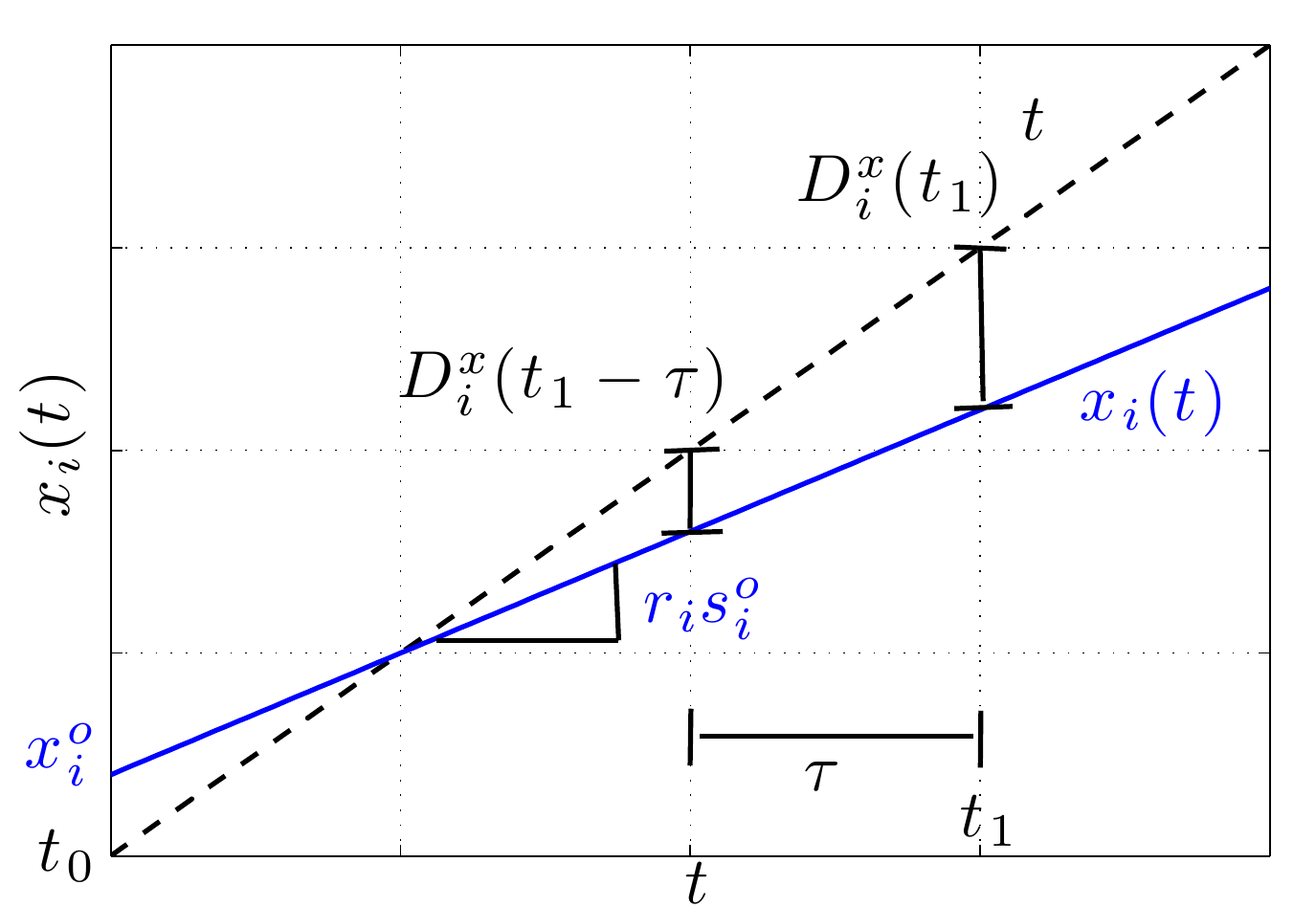}
		\caption{ Offset and relative skew measurements }\label{fig:xt3}
        \end{subfigure}%
        \vspace{-.25cm}
        \caption{Time estimation and relative measurements}\label{fig:fig1}
\end{figure}

Thus, synchronizing computer clocks implies correcting $x_i(t)$ in order to match $t$, i.e. enforcing $x_i(t)=t$.
There are two difficulties on this estimation process. Firstly, the initial time $t_0$ in which the counter starts its unknown. Secondly,
the counter updating period $\delta_i$ ($\delta_i\approx\delta^0$) is usually unknown with enough precision and therefore presents a skew  $r_i=\frac{x_i(t)-x_i(t_0)}{t-t_0}=\frac{\delta^0}{\delta_i}$.
This is illustrated in Figure \ref{fig:xt} where $x_i(t)$ not only increases at a different rate than $t$, but also starts from a value different from $t_0$, represented by  $x_i^o$.

In practice, $c_i(t)$ can be approximated by a real \changed{value} since the time between increments is extremely small ($0.416$ns) and the maximum count register value so large ($2^{64}-1$) that it would take more than $200$ years to reach.
Therefore, $x_i(t)$ can be described by the linear map of the global reference $t$, i.e.
\begin{equation}\label{eq:linear_map}
x_i(t) = r_is_i^o(t-t_0) +x_i^o,
\end{equation}
where $s_i^o$ is an additional skew correction implemented to compensate the skew.  Equation \eqref{eq:linear_map} also shows that if one can set $s_i^o=1/r_i$ and $x_i^o=t_0$, then we obtain a perfectly synchronized clock with $x_i(t)=t$.


\begin{figure}[htp]
       \begin{subfigure}[b]{0.53\columnwidth}
               \centering
               \includegraphics[width=\columnwidth,height = 3cm]{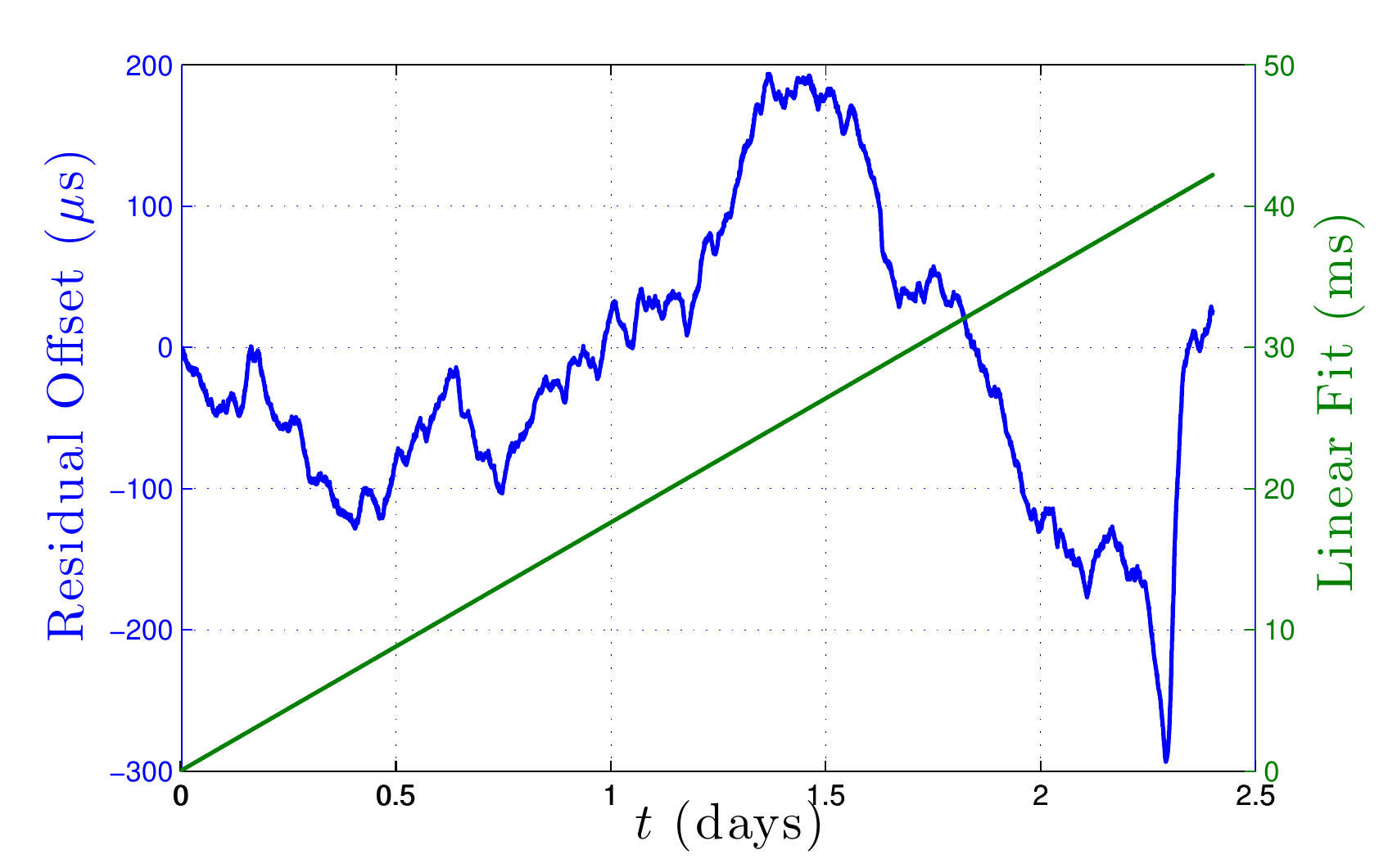}
               \caption{\changed{Variation of the offset between two TSC counters changes on skew ($r_i$): The right y-axis represents the mean offset change in ms and the left y-axis the residual offset (offset minus the mean) in $\mu$s}}\label{fig:tsc_drift}
        \end{subfigure}\;\;
        \begin{subfigure}[b]{0.43\columnwidth}
               \centering
		\includegraphics[width=\columnwidth,height =3cm]{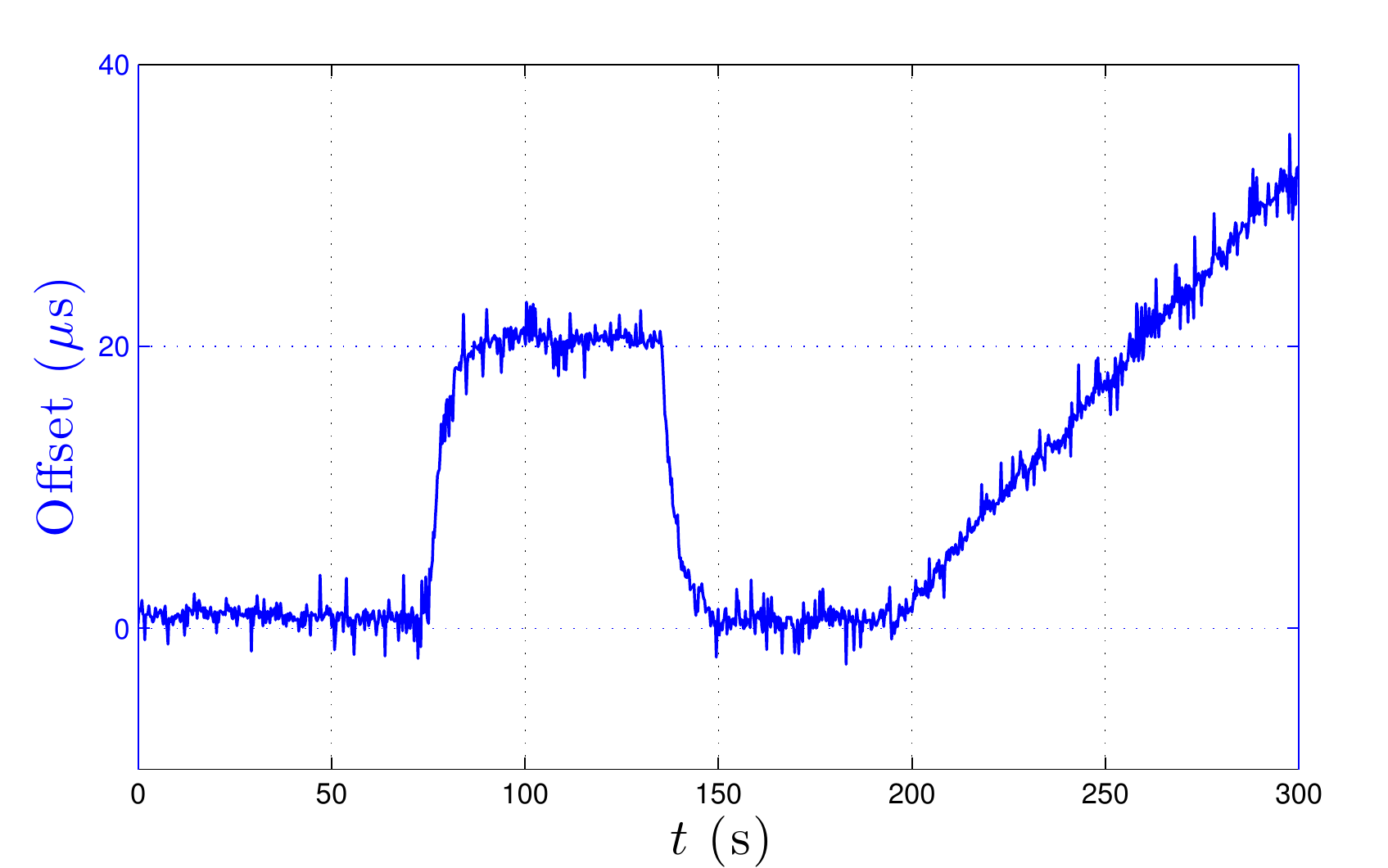}
		\caption{ Example of skew and offset corrections on linux time: First a $20\mu$s offset is added and subtracted and then a skew of $0.3$ppm is introduced}\label{fig:adjtimex}
        \end{subfigure}%
        \vspace{-.25cm}
        \caption{Comparison between two TSC counters, and skew and offset corrections using {\it adjtimex()}}\label{fig:tsc_and_adjtimex}
\end{figure}

The main problem is that not only neither $t_0$ nor $r_i$ can be explicitly estimated, but also $r_i$ varies with time as shown in Figure \ref{fig:tsc_drift}.
Thus, current protocols periodically update $s_i^o$ and $x_i^o$ in order to keep track of the changes of $r_i$. These updates are made using the {\it offset} between the current estimate $x_i(t)$ and the global time $t$, i.e. $D_i^x(t)=t-x_i(t)$, and the {\it relative frequency error} that is computed using two offset measurements separated by $\tau$ seconds, i.e.
\begin{equation}\label{eq:ferr}
f^{err}_i(t):=\frac{D_i^x(t)-D_i^x(t-\tau)}{x_i(t)-x_i(t-\tau)} = \frac{1-r_is_i^o}{r_is_i^o}.
\end{equation}
Figure \ref{fig:xt3} provides an illustration of these measurements.
\changed{In most protocols (see e.g. ~\cite{mills_network_2010,froehlich_achieving_2008,Veitch:2009td}) \eqref{eq:ferr} goes through an additional filtering process to reduce the estimation noise. Here we will use $f^{err}_i(t_k)$ to denote either the measurement obtained using \eqref{eq:ferr} or a filtered version of it.}



%

To understand the differences between  current protocols, we first rewrite the evolution of $x_i(t)$ based only on the time instants $t_k$ in which the clock corrections are performed. We allow the skew correction $s_i^o$ to vary over time, i.e. $s_i(t_k)$, and write $x_i(t_{k+1})$ as a function of $x_i(t_{k})$. Thus, we obtain
\begin{subequations}
\label{eq:double_int}
\begin{align}
\centering
x_i(t_{k+1}) &= x_i(t_k) +  \tau r_is_i(t_k)  + u_i^x(t_k)\label{eq:double_int_ia}\\
s_i(t_{k+1}) &=s_i(t_k) +  u_i^s(t_k)\label{eq:double_int_ib}
\end{align}
\end{subequations}
where $\tau = t_{k+1}-t_k$ is the time elapsed between adaptations; also known as poll interval~\cite{mills_network_2010}.
The values $u_i^x(t_k)$ and  $u_i^s(t_k)$ represent two different types of corrections that a given protocol chooses to do at time $t_k$ and are usually implemented within the interval $(t_k,t_{k+1})$. $u_i^x(t_k)$ is usually referred to as {\it offset correction } and $u_i^s(t_k)$ as  {\it skew correction}.\footnote{These corrections can be implemented in Linux OS using the {\it adjtimex()} interface to update the system clock or by maintaining a virtual version of $x_i(t)$ and directly applying the corrections to it, as in IBM CCT~\cite{froehlich_achieving_2008} and RADclock~\cite{Veitch:2009td}. The latter gives more control on how the corrections are implemented since it does not depend on kernel's routines.}
See Figure \ref{fig:adjtimex} for an illustration of their effect on the linux time.
%
\begin{remark}
\changed{One of the implicit assumptions of the model \eqref{eq:double_int} is that we require every server to update their clocks simultaneously at time instances $\{t_k\}$. This may seem unrealistic since its implementation would require sharing a common time reference which is the whole purpose of the algorithm. However, the analysis presented in Section \ref{sec:analysis} can be extended for pseudo-synchronous implementations as proposed in \cite{Simeone:2008cc} where each node measures the offset with their neighbors and updates whenever $x_i(t)=kT$.}
\end{remark}

\changed{We now proceed to summarize the different types of adaptations implemented by current protocols. To simplify the comparison, we assume that each server can connect directly to the source of UTC time ($t$). This assumption will be dropped in Section \ref{sec:algorithm} after we describe our solution.}
 The main differences between current protocols lies on whether they use offset corrections, skew corrections, or both, and whether they update using offset values $D_i^x(t_k)$, relative frequency errors $f_i^{err}(t_{k})$, or both.

\subsection{ Offset corrections}\label{sssec:only_offset}
These corrections consist in keeping the skew fixed and periodically introducing time changes of size  $u_i^x(t_k) = \kp_1 D_i^x(t_k)$ or $u_i^x(t_k) = \kp_1 D_i^x(t_k) + \kp_2 f_i^{err}(t_k)$
where $\kp_1,\kp_2>0$.
They are used by NTPv3~\cite{mills_network_1992}  and NTPv4~\cite{mills_network_2010} respectively under ordinary conditions.

\begin{figure}
       \begin{subfigure}[b]{0.49\columnwidth}
               \centering
               \includegraphics[width=.9\columnwidth,height = .5\columnwidth]{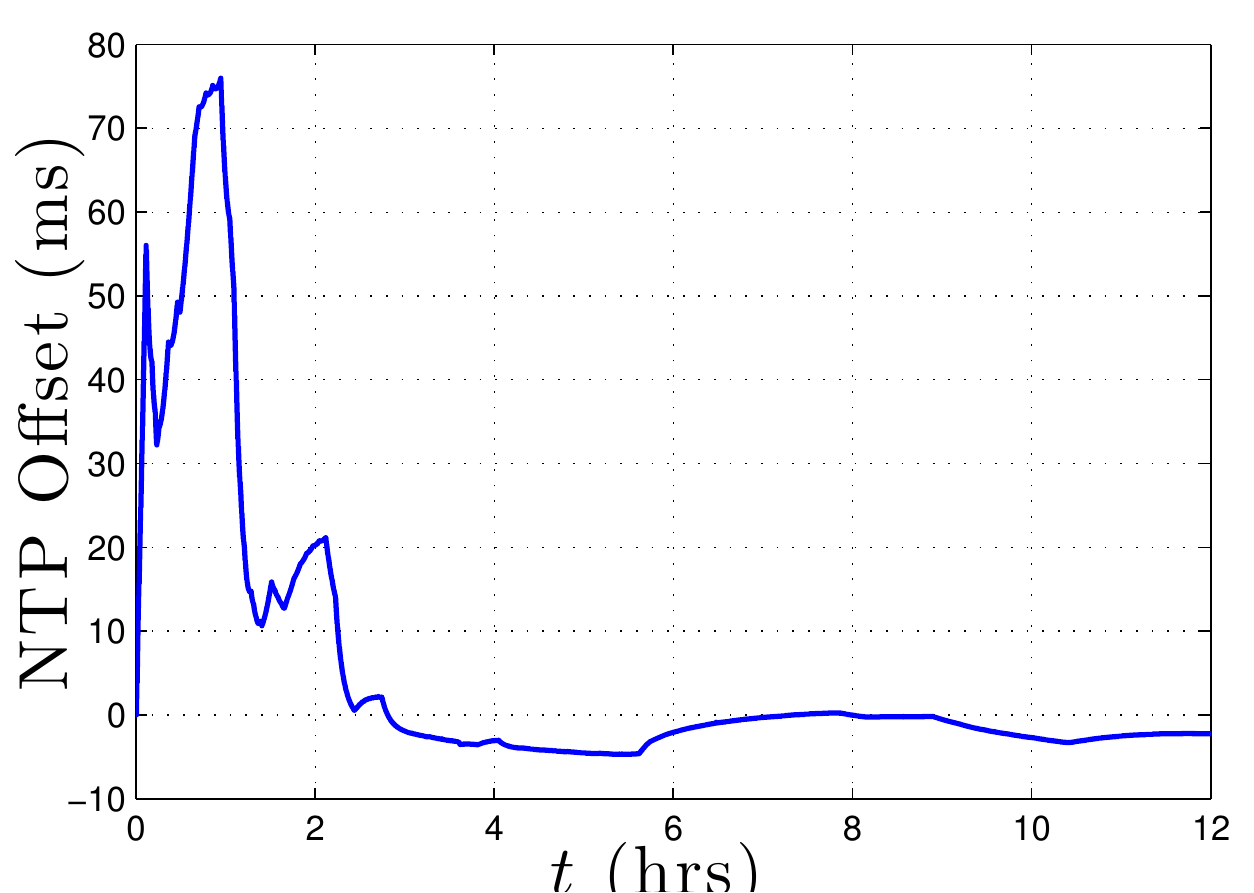}
               \caption{NTP initialization period}\label{fig:ntp_init}
        \end{subfigure}
        \begin{subfigure}[b]{0.49\columnwidth}
               \centering
		\includegraphics[width=.9\columnwidth,height = .5\columnwidth]{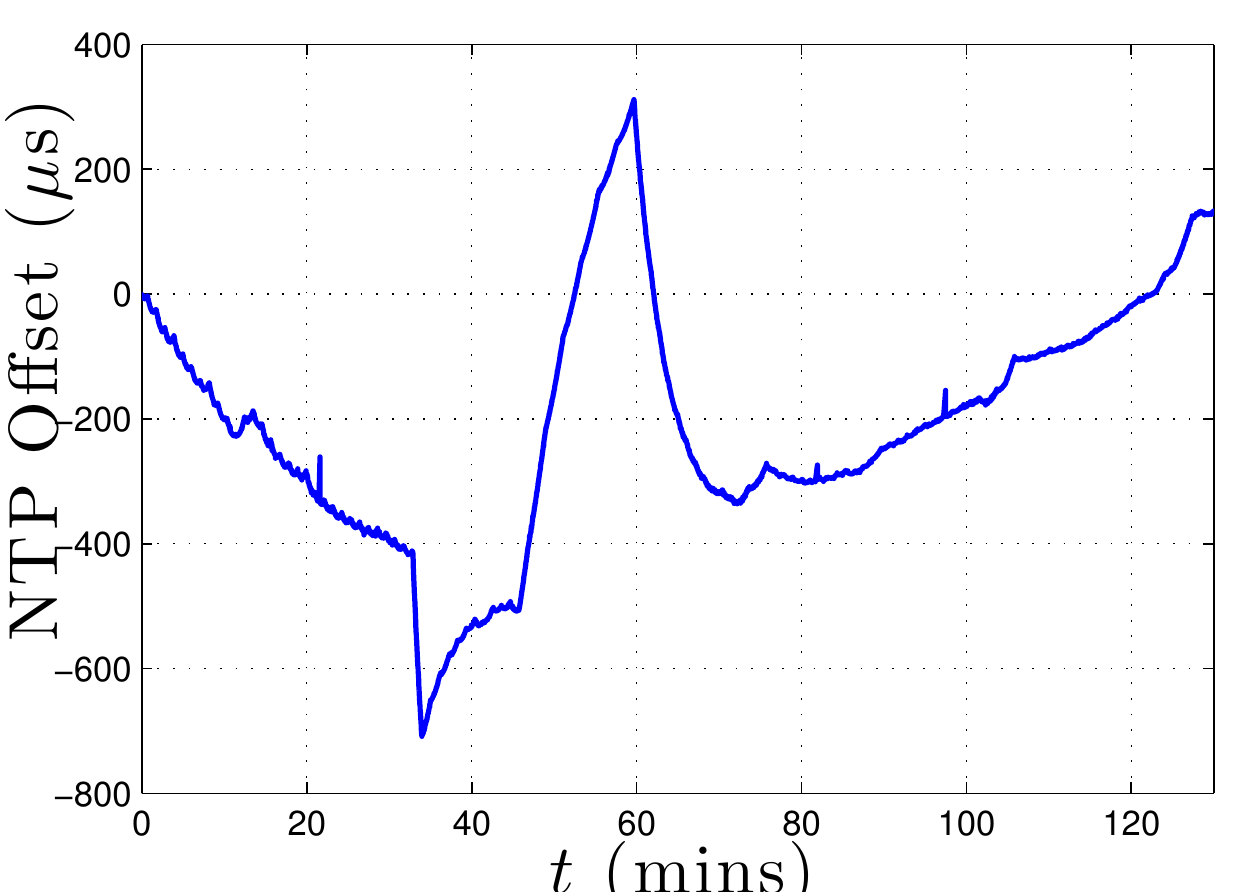}
               \caption{NTP in normal regime}\label{fig:ntp_normal}
        \end{subfigure}%
        \vspace{-.25cm}
        \caption{ Variations of NTP time using TSC as reference }\label{fig:ntp}
\end{figure}

\begin{figure}
       \begin{subfigure}[b]{0.33\columnwidth}
               \centering
               \includegraphics[width=\columnwidth]{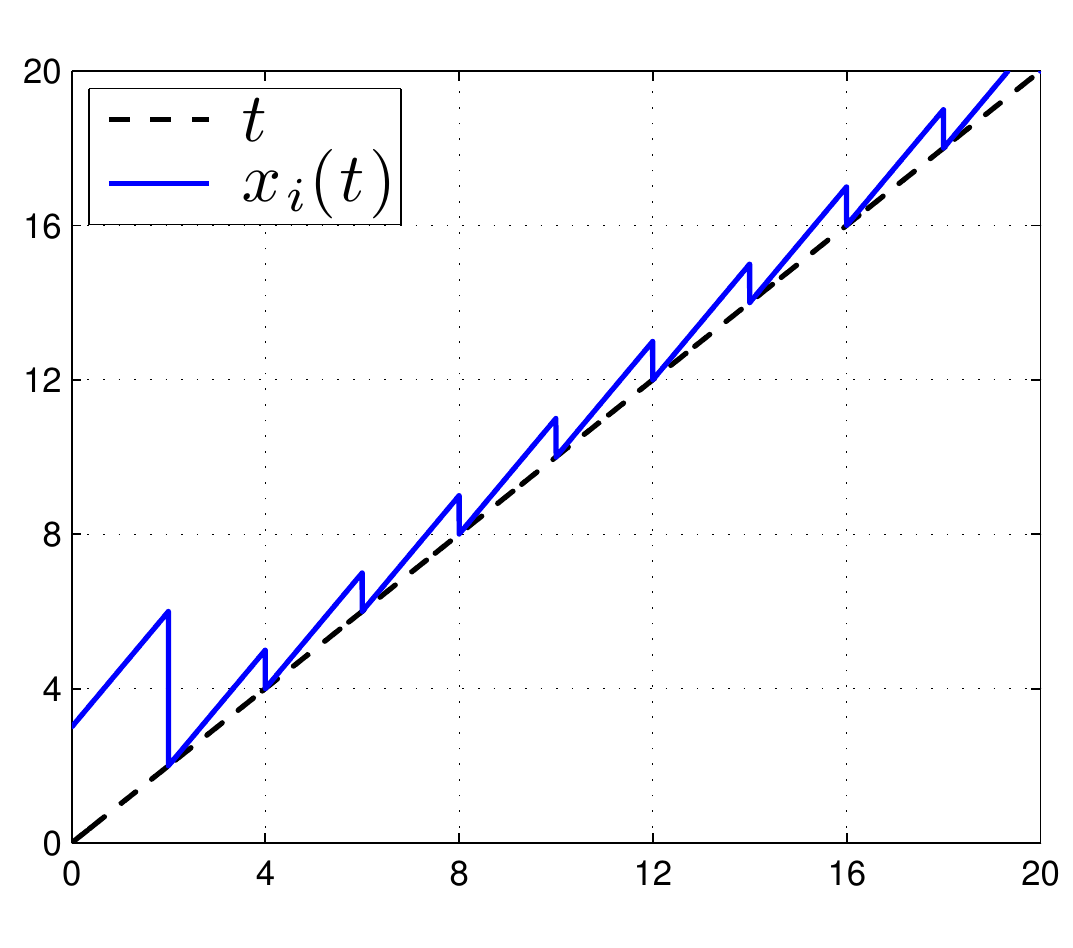}
               \caption{\text{Offset corrections}\quad\text{        }\quad\quad}\label{fig:alg-offset}
        \end{subfigure}
        \begin{subfigure}[b]{0.33\columnwidth}
               	\centering
               	\includegraphics[width=\columnwidth]{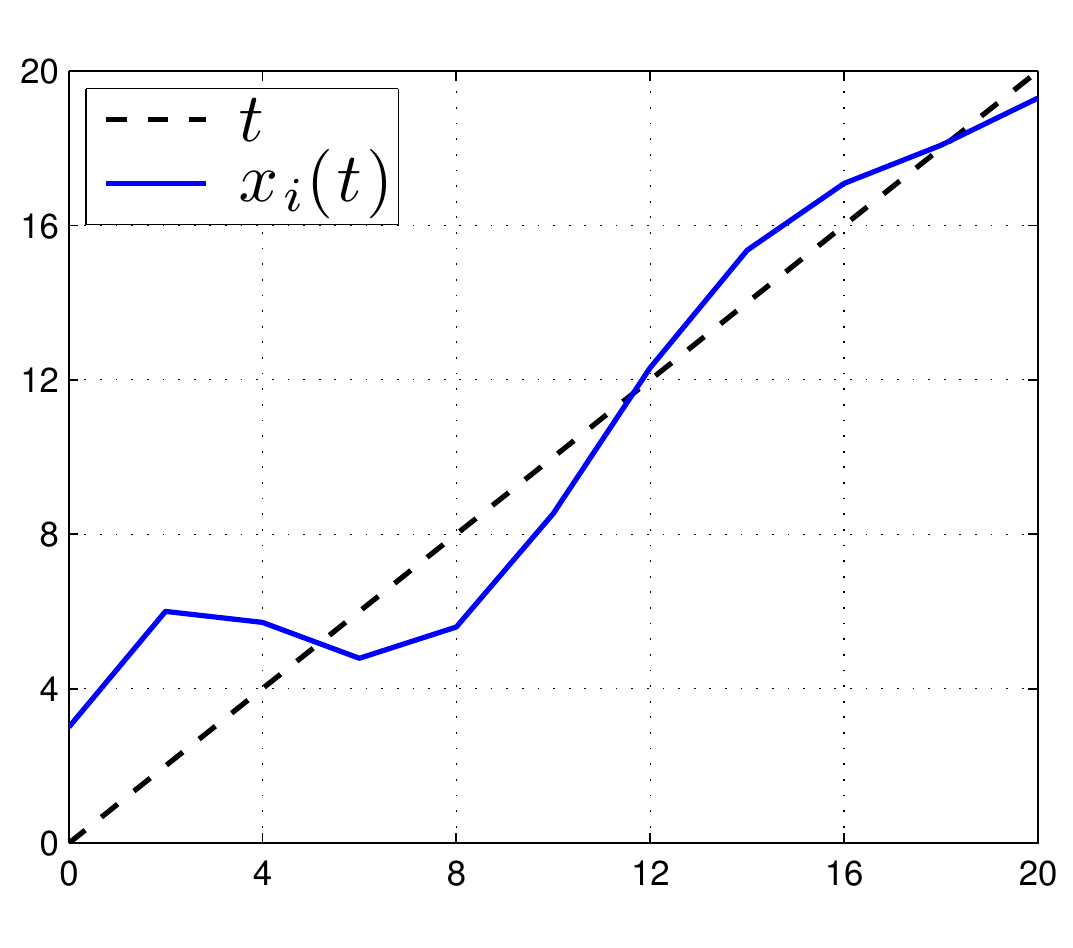}
				\caption{\text{Skew corrections}\quad\text{           }\quad\quad}\label{fig:alg-skew}
        \end{subfigure}%
        \begin{subfigure}[b]{0.33\columnwidth}
               	\centering
               	\includegraphics[width=\columnwidth]{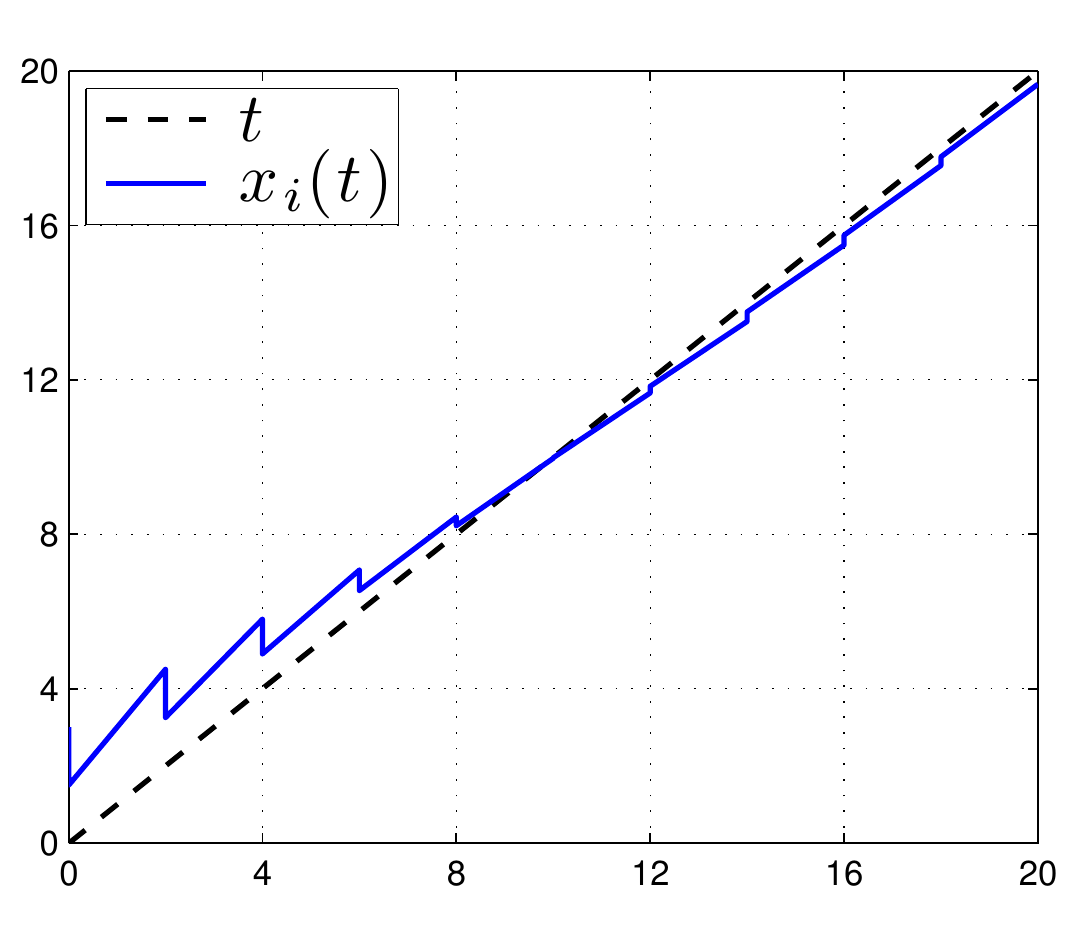}
				\caption{Offset and skew corrections}\label{fig:alg-offsetskew}
        \end{subfigure}%
        \vspace{-.25cm}
        \caption{Current Protocols Adaptation}\label{fig:alg}
\end{figure}

These protocols have in general a slow initialization period as shown in Figure \ref{fig:ntp_init}. This is because the algorithm must first obtain a very  accurate estimate of the initial frequency error $f_i^{err}(t_0)$. 
Furthermore, these updates usually generate non-smooth time evolutions as in Figures \ref{fig:ntp_normal}  and \ref{fig:alg-offset}, and should be done carefully since they might introduce backward jumps ($x_i(t_{k+1})<x_i(t_k)$), which can be problematic for some applications.

\subsection{Skew corrections}\label{sssec:only_skew}
Another alternative that avoids using steep changes in time is proposed by the IBM CCT solution~\cite{froehlich_achieving_2008}. This alternative does not introduce any offset correction, i.e. $u_i^x(t_k) = 0$, and updates the skew $s_i(t_k)$ by
$u_i^s(t_k) = \kappa_1 D_i^x(t_k) + \kappa_2 f_i^{err}(t_k)$.

The behavior of this algorithm is shown in Figure \ref{fig:alg-skew}. In \cite{xie_consensus_2011} it was shown for a slightly modified version of it 
(used $r_is_i(t_k)f_i^{err}(t_k)$ instead of $f_i^{err}(t_k)$) the algorithm can achieve synchronization for very diverse network architectures.

However, the estimation of $f_i^{err}(t_k)$ is nontrivial as it is constantly changing with subsequent updates of $s_i(t_k)$ and it usually involves sophisticated computations~\cite{zhang_clock_2002,Kim:2012db}.



\subsection{Skew and offset corrections}

This type of corrections allow dependence on only offset information $D_i^x(t_k)$ as input to $u_i^x(t_k)$ and $u_i^s(t_k)$. For instance, in \cite{Carli:2014gd} the update
$
u_i^x(t_k) = \kappa_1D_i^x(t_k) \text{ and } 
u_i^s(t_k) = \kappa_2D_i^x(t_k)
$
was proposed.
This option allows the system to achieve synchronization without any skew estimation. But the cost of achieving it, is introducing  offset corrections in $x_i(t)$ as shown in Figure \ref{fig:alg-offsetskew}. Therefore, it suffers from the same problems discussed in \ref{sssec:only_offset}.

\changed{Another alternative that falls in into this category is the RADclock~\cite{Veitch:2009td}. In this solution the offset correction $u_i^x(t_k)$ is an exponential average of the past offsets and the skew compensation $u_i^s(t_k)$ is a filtered version of $f_i^{err}(t_k)$. The exponential average of offsets and filter stage in $f_i^{err}(t_k)$ allows this solution to mitigate the jumps and become more robust to jitter. However, it does not necessarily prevent backward jumps unless the offset corrections are smaller than the precision of the clock.}

\section{Continuous Skewless Synchronization}\label{sec:algorithm}

We now present an algorithm that overcomes the limitations of the solutions described in Section \ref{sec:comp_clocks}. In other words, our solution  has the following two properties:
\begin{enumerate}
\item Continuity: The protocol does not introduce steep changes on the time value, i.e. $u_i^x(t_k)\equiv 0$.
\item Skew independence: The protocol does not use skew information $f_i^{err}(t_k)$ as input.
\end{enumerate}
A solution with these properties will therefore prevent unnecessary offset corrections that produce jitter and will be more robust to noise by avoiding skew estimation.
After describing and motivating our algorithm, we show how the updating rule can be  implemented in the context of a network environment.

\begin{figure}
\centering
\includegraphics[height=.3\columnwidth,width=.85\columnwidth]{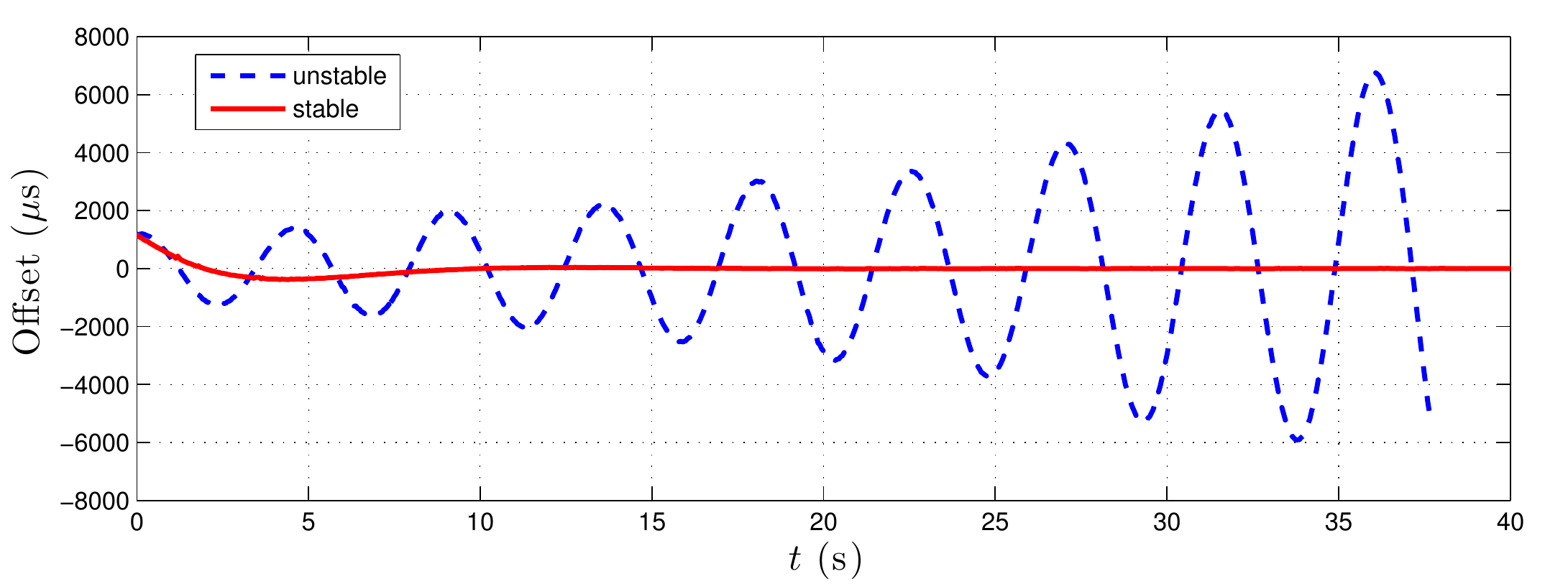}
\caption{Unstable clock steering using only offset information \eqref{eq:double_int_pure} and stable clock steering based on exponential average compensation~\eqref{eq:algorithm} }\label{fig:unstable}
\end{figure}

The motivation behind the proposed solution comes from trying to compensate the problem that arises when one tries to naively impose properties 1) and 2), i.e. using
\begin{align} \label{eq:double_int_pure}
u_i^x(t_k) = 0 \quad \text{ and } \quad u_i^s(t_k) = \kp_1 D_i^x(t_k).
\end{align}
Figure \ref{fig:unstable} shows that this type of clock corrections is unstable; the offset $D_i^x(t_k)$ of the slave clock oscillates with an exponentially increasing amplitude.

The oscillations in Figure \ref{fig:unstable} arise due to the fundamental limitations of using offset to update frequency.
This is better seen in the continuous time version of the system \eqref{eq:double_int} with \eqref{eq:double_int_pure}, i.e.
\[
\dot x_i(t) = r_i s_i(t)\text{ and } \dot s_i(t) = \kp_1 D_i^x(t)
\]
where $\dot x(t) =\frac{d}{dt}x(t)$.
If we consider the offset $D_i^x=t-x_i(t)$ as the system state, then we have
\[
\dot D_i^x =1-r_is_i\text{ and }
\ddot D_i^x=-\kp_1 r_iD_i^x,
\]
with $\ddot x(t)=\frac{d^2}{dt^2} x(t)$.

This is analogous to a spring mass system without friction. Thus, it has two purely imaginary eigenvalues that generate sustained oscillations; see \cite{Mallada:2008jb,Mallada:2011iu} for similar examples.\footnote{In the discrete time system the oscillations increase in amplitude since there is a delay between the time the offset is measured $t_k$ and the time the update is made $t_{k+1}$ which makes the system unstable.}
One way to damp these oscillations in the spring-mass case is by adding {\it friction}. This implies  adding a term that includes a frequency mismatch $f_i^{err}(t)$ in our system, which is
equivalent to the protocols of Section \ref{sssec:only_skew}, and therefore
undesired.

However, there are other ways to damp these oscillations using passivity-based techniques from control theory~\cite{
ren2008distributed}. The basic idea is to introduce an additional state $y_i$ that generates the desired {\it friction} to damp the oscillations.

Inspired by \cite{ren2008distributed}, we consider the exponentially weighted moving average of the offset
\begin{equation}\label{eq:y}
y_i(t_{k+1}) = pD_i^x(t_k) + (1-p)y_i(t_{k}).
\end{equation}
and update $x_i(t_k)$ and $s_i(t_k)$ using:
%
\begin{align}
u_i^x(t_{k}) = 0 \quad\text{ and }\quad u_i^s(t_{k}) &= \kappa_1 D^x(t_k) -\kappa_2 y(t_{k}).\label{eq:algorithm}
\end{align}
Figure \ref{fig:unstable} shows how the proposed strategy is able to compensate the oscillations without needing to estimate the value of $f_i^{err}(t_k)$.
The stability of the algorithm will depend on how $\kp_1$, $\kp_2$ and $p$ are chosen. A detailed specification of these values is given in Section \ref{ssec:param_vals}.

Finally, since we  are interested in studying the effect of timing loops, we move away from the client-server configuration implicitly assumed in Section \ref{sec:comp_clocks} and allow mutual or cyclic interactions among nodes.
The interactions between different nodes is described by a  graph $G(V,E)$, where $V$ represents the set of $n$ nodes ($i\in V$) and $E$ the set of {\it directed} edges $ij$; $ij\in E$ means node $i$ can measure its offset with respect to $j$, $D_{ij}^x(t_k) = x_j(t_k)-x_i(t_k)$.

Within this context, a natural extension of  \eqref{eq:y}-\eqref{eq:algorithm} is to substitute $D_i^x(t_k)$ with the weighted average of $i$'s neighbors offsets.
Thus, we propose the following algorithm to update the clocks in the network.

\vspace{.115cm}
\noindent{\bf Algorithm 1 (Alg1):}
For each computer node $i$ in the network, perform the following actions:\noindent
\begin{itemize}[-]
\item  Compute the time offsets ($D^x_{ij}(t_k)$) from $i$ to every neighbor $j$ at time $t_k$.
\item  Update the skew $s_i(t_{k+1})$ and the moving average $y_i(t_{k+1})$ at time $t_{k+1}$ according to:
\begin{subequations}
 \label{eq:system}
 \begin{align}
x_i(t_{k+1}) =& x_i(t_{k}) +  \tau_{k}r_is_i(t_{k})\label{eq:system_xi}\\
s_i(t_{k+1}) =&s_i(t_{k}) +  \kp_1\sum_{j\in \mathcal N_i} \alpha_{ij}D^x_{ij}(t_k)
- \kappa_2y_i(t_{k}) \label{eq:system_si}\\
y_i(t_{k+1}) =& p\sum_{j\in \mathcal N_i}\alpha_{ij}D^x_{ij}(t_k) + (1-p)y_i(t_{k})\label{eq:system_yi}
 \end{align}
\end{subequations}
where $\mathcal N_i$ represents the set of neighbors of $i$ and the weights $\alpha_{ij}$ are positive.
\end{itemize}

\vspace{.115cm}
\changed{Equation \eqref{eq:system} can be interpreted as a discrete-time {\it second-order consensus} algorithm with an additional smoothing in which, besides using position information (time estimates $x_i(t_k)$), we use a smoothed version of the position errors ($y_i(t_k)$) to control speed ($s_i(t_k)$). Consensus algorithms have been a subject of intensive research since the seminar work of Jadbabaie et al.~\cite{Jadbabaie:2003bu}, see e.g. \cite{ren2008distributed} and references therein. In particular, application of consensus ideas to computer clock synchronization can be found in~\cite{Carli:2014gd} (second order consensus) and \cite{Simeone:2008cc} (first order consensus). Thus, the analysis presented in this paper also contributes to this rich literature by characterizing convergence of discrete-time consensus algorithms.}

When using our algorithm, many servers can affect the final frequency of the system. Thus, when the system synchronizes, we have
\begin{equation}\label{eq:synchronization}
\changed{x_i(t_k) \rightarrow x^{\text{ref}}(t_k):= r^*(t_k-t_0) + x^*\quad i\in V.}
\end{equation}
$r^*$ and $x^*$ are possibly different from their ideal values $1$ and $t_0$. Their final values depend on the initial condition of all different clocks as well as the topology, which we assume to be a connected graph in this paper.





\subsubsection*{Differences with RADclock}
\changed{Although \eqref{eq:system} seems to be similar to RADclock~\cite{Veitch:2009td}, there are some key differences that affect their behavior.\\
\noindent 1) Even though both solutions used an exponentially weighted offset estimate, our filtering  \eqref{eq:system_yi} does not depend on the estimated offset error as in~\cite{Veitch:2009td}. Moreover, while RADclock uses it to make offset corrections (changing $u_i^x(t_k)$), we use our weighted offset measurement $y_i(t)$ to make skew correction (changing $u_i^s(t_k)$). Therefore, neither the measurement itself nor its use are the same.\\
\noindent 2) RADclock explicitly uses offset measurements to introduce correction on the offset ($u_i^x(t_k)$) and an estimation of the skew to compensate it  ($u_i^s(t_k)$). Our algorithm only compensates the skew by using the last measured offsets $D_{ij}^x(t_k)$ and our filtered offset measurement $y_i(t_k)$. Thus, we have neither explicit estimation of the skew nor explicit compensation of the offset, which makes synchronization rather unintuitive. 
}

\subsubsection*{Notation}
We use $\0_{m\times n}$ ($\1_{m\times n}$) to denote the matrices of all zeros (ones) within $\mathds R^{m\times n}$ and  $\0_{n}$ ($\1_{n}$) to denote the column vectors of appropriate dimensions. $I_n\in\mathds R^{n\times n}$ represents the identity matrix.  Given a matrix $A\in\mathds R^{n\times n}$ with Jordan normal form $A=PJP^{-1}$, let $n_A\leq n$ denote the total number of Jordan blocks $J_l$ with $l\in\mathcal I(A):=\{1,...,n_A\}=|\mu_1(A)|$. 
We use $\mu_l(A)$, $l\in\{1,\dots,n\}$ or just $\mu(A)$ to denote the eigenvalues of $A$, and order them decreasingly $|\mu_{1}(A)|\geq\dots\geq |\mu_{n}(A)|$. The function $\rho(A)$ is the spectral radius of $A$ or equivalently the largest absolute value of its eigenvalues $\rho(A)=\max_{l\in\mathcal I(A)}|\mu_l(A)|$.
 Finally, $A^T$ is the transpose of $A$, $A_{ij}$ is the element of the $i$th row and $j$th column of $A$ and $a_i$ is the $i$th element of the column vector $a$, i.e. $a=[a_i]^T$.

\section{Convergence Analysis}\label{sec:analysis}

We now analyze the asymptotic behavior of system \eqref{eq:system} and provide a necessary and sufficient condition on the parameter values that guarantee its convergence to \eqref{eq:synchronization}. 
\changed{ Throughout this section we shall assume that the internal skew $r_i$ of each clock is constant and that the offset measurements $D_{ij}^x(t_k)$ can be obtained without incurring in any error. These assumptions will be relaxed in Section \ref{sec:noise}.
}

\changed{
The key insight of our analysis comes from decomposing the system \eqref{eq:system_z} into two complementary systems that keep track of two different physical properties. In particular, we will use the scalars
\begin{gather}\label{eq:average}
\begin{aligned}
\tilde x(t_k) := \gamma&\sum_{i=1}^{n} \frac{\xi_i}{r_i} x_i(t_k),\quad \tilde s(t_k) :=\gamma\sum_{i=1}^{n} \xi_i s_i(t_k)\\
&\text{ and } \quad\tilde y(t_k) :=\gamma\sum_{i=1}^{n} \xi_i y_i(t_k)
\end{aligned}
\end{gather}
to track the average behavior of the system and 
\begin{gather}\label{eq:deviations}
\begin{aligned}
\delta x_i(t_k):= x_i&(t_k) - \tilde x(t_k),\; \delta s_i(t_k):= s_i(t_k) - \frac{\tilde s(t_k)}{r_i} \\
 &\text{ and }\quad  \delta y_i(t):=y_i(t)- \frac{\tilde y(t_k)}{r_i}
 \end{aligned}
\end{gather}
to track how each individual clock deviates from the collective mean.
}
\changed{
Here, $\xi=[\xi_i]^T$ is the normalized ($\sum_i \xi_i = 1$) left eigenvector of the zero eigenvalue of the Laplacian matrix $L\in \mathds R^n$ associated with $G(V,E)$, i.e.
\begin{align*}
L_{ii} =\alpha_{ii}:= \sum_{j\in\mathcal N_i} \alpha_{ij} \text{ and } L_{ij} &=
\begin{dcases*}
 -\alpha_{ij}                        & if $ij \in E$,\\
0                      & otherwise.
\end{dcases*}
\end{align*}
and $\gamma$ is the $\xi_i$-weighted harmonic mean of $r_i$, i.e. $\frac{1}{\gamma} = \1_n^T R^{-1}\xi = \sum_{i=1}^{n} \frac{\xi_i}{r_i}.$ While in general $\xi$ may not be unique, it becomes unique when $G(E,V)$ is connected.}

\changed{It will be more convenient sometimes to use a vector form representation of \eqref{eq:system} given by}
\begin{equation}\label{eq:system_z}
z_{k+1} = A z_k
\end{equation}
where $z_k:=[ x(t_{k})^T s(t_{k})^T y(t_{k})^T]^T\in\mathds R^{3n}$,
\begin{equation*}
A:=
\left[
\begin{array}{ccc}
I_n & \tau R &\0_{n\times n} \\
-\kappa_1L & I_n & -\kappa_2I_n \\
p(-L) & \0_{n\times n} & (1-p) I_n
\end{array}
\right]\in \mathds{R}^{3n\times 3n},
\end{equation*}
$R\in\mathds R^{n\times n}$ is the diagonal matrix with elements $r_i$.
\changed{ Similarly, we can express the evolution of $\tilde z_k:=[\tilde x(t_k)\; \tilde s(t_k)\; \tilde y(t_k)]^T$ and $\delta z_k:=[\delta x(t_k)^T\; \delta s(t_k)^T\;
\delta y(t_k)^T]^T$ using}
\begin{equation}\label{eq:system_split}
\changed{\delta z_{k+1} = \hat A\delta z_k\quad \text{ and }\quad \tilde z_{k+1} = \tilde A\tilde z_k},
\end{equation}
\changed{where $\hat A := NA$, $N:=\blockdiag(N_1,N_2,N_2)$, $N_1:=I_n - \gamma\1_n\xi^TR^{-1}$, $N_2:=I_n-\gamma R^{-1}\1_n\xi^T$ and }
\begin{equation}
\changed{\tilde A:=\left[
\begin{array}{ccc}
1 & \tau & 0\\
0 & 1 &-\kp_2\\
0 & 0 &1-p
\end{array}
\right].}
\end{equation}

The convergence analysis of this section is done in two stages. First, we provide necessary and sufficient conditions for synchronization in terms of the eigenvalues of $A$ (Section \ref{ssec:asymptotic_behavior}) and then use Hermite-Biehler Theorem \cite{bhattacharyya_robust_1995} to relate these eigenvalues with the parameter values  that can be directly used in practice (Section \ref{ssec:param_vals}). All the proof details are included in the appendix for interested readers.

\subsection{Asymptotic Behavior}\label{ssec:asymptotic_behavior}
We start  by studying the asymptotic behavior of \eqref{eq:system_z}. That is, we are interested in finding under what conditions the series of elements $\{x_i(t_k)\}$ converge to \eqref{eq:synchronization} as $t_k$ goes to infinity.

\changed{We will show that we can study \eqref{eq:system_z} by looking at the evolution of \eqref{eq:system_split}. In particular we will show that \eqref{eq:synchronization} is equivalent to}
\begin{subequations}\label{eq:synchronization_split}
\begin{gather}
\changed{\delta x(t_k) \rightarrow \0_{n}, \quad \delta s(t_k) \rightarrow \0_{n},\quad \delta y(t_k) \rightarrow \0_{n}, \label{eq:deltaz_sync}}\\
\changed{\tilde x(t_{k}) \rightarrow x^{\text{ref}}(t_k), \quad \tilde s(t_k) \rightarrow r^*\text{ and }\;\tilde y(t_k)\rightarrow 0. \label{eq:tildez_sync}}
\end{gather}
\end{subequations}


Consider the Jordan normal form \cite{horn_matrix_1990} of 
\begin{equation}\label{eq:jordan_form}
A = PJP^{-1}
    := \left[ \zeta_1 \quad ... \quad\zeta_{3n}\right] 
J
\left[ \eta_1 \quad ... \quad\eta_{3n}\right]^T
\end{equation}
 where $J=\blockdiag(J_l)_{l\in \mathcal I(A)}$, $\zeta_i$ and $\eta_i$ are the right and left generalized eigenvectors of $A$ such that
\begin{align*}
\zeta_i^T\eta_j =
\begin{dcases*}
1 & if $j=i$,\\
0 & otherwise.
\end{dcases*}
\end{align*}
The following lemmas will allow us to connect the behavior of \eqref{eq:system_split} with \eqref{eq:system_z}.

\begin{lemma}[Eigenvalues of $A$ and Multiplicity of $\mu(A) = 1$] \label{lem:multiplicity}
$A$ has an eigenvalue $\mu(A)=1$ with multiplicity $2$ if and only if the graph $G(V,E)$ is connected, $\kappa_1\neq\kappa_2$ and $p>0$.

Furthermore, $\mu_l(A)$ are the roots of
\changed{
\begin{equation}\label{eq:g_l}
g_l(\lambda)=(\lambda-1)^2(\lambda-1+p) +[(\lambda-1)\kappa_1 +p(\kappa_1-\kappa_2)]\nu_l
\end{equation}
}
where $\nu_l = \mu_l(\tau LR)$ and satisfies
\begin{equation}\label{eq:nu_condition}
\nu_n=0<|\nu_l|\text{ for }l\in\{1,\dots,n-1\}.
\end{equation}
\end{lemma}


\begin{lemma}[Jordan Chains Properties]\label{lem:jordan_chains}
Under the conditions of Lemma \ref{lem:multiplicity} the right and left Jordan chains, $(\zeta_1,\zeta_2)$ and $(\eta_2,\eta_1)$ respectively, associated with $\mu(A)=1$ and the eigenvectors $\zeta_3$ and $\eta_3$ associated with $\mu(A)=1-p$ are given by
\begin{equation}\label{eq:zeta}
[\zeta_1\;\zeta_2\;\zeta_3]= \left[\begin{array}{ccc}
\1_n & \1_n                                & -\frac{\tau\kp_2}{p^2}\1_n \\
\0_n & \frac{(R^{-1}\1_n)}{\tau} & \frac{\kp_2}{p}R^{-1}\1_n \\
\0_n & \0_n                                & R^{-1}\1_n
\end{array}
\right] \text{ and }
\end{equation}
\changed{\begin{equation}\label{eq:eta}
[\eta_1\;\eta_2\;\eta_3]=\gamma \left[\begin{array}{ccc}
 R^{-1}\xi 																&	\0_n							& 	\0_n		\\
 -\tau\xi 																&	\tau\xi								&	\0_n		\\
 \tau\kp_2(\frac{1}{p}+\frac{1}{p^2})\xi 	& -\tau\frac{\kp_2}{p}\xi 		&	\xi
\end{array}
\right].
\end{equation}}
\changed{Moreover, given $\zeta_l=[x_l^T\; s_l^T\; y_l^T]^T$, with $l>3$, then the following conditions must be satisfied
\begin{equation}\label{eq:zeta_l_condition}
\1_n\xi^TR^{-1}x_l = \0_n,\;\;  \1_n\xi^Ts_l=\0_n\;\text{ and }\;  \1_n\xi^Ty_l=\0_n.
\end{equation}}
\end{lemma}


The proof of Lemmas \ref{lem:multiplicity} and \ref{lem:jordan_chains} can be found in Appendices \ref{app:lemma1} and \ref{app:lemma2}.

\changed{
Lemmas \ref{lem:multiplicity} and \ref{lem:jordan_chains} also provide further information of the structure of $J$ in \eqref{eq:jordan_form}. That is,  $J=\blockdiag(\hat J_1, \hat J_2)$
where 
\begin{equation}\label{eq:hatJ_1}
\hat J_1 := \left[   \begin{array}{crc}
1 & 1 & 0\\
0 & 1 & 0\\
0 & \;\; 0 & 1-p
\end{array}
\right].
\end{equation}
Moreover, direct application of \eqref{eq:zeta_l_condition} shows that $N=P[\blockdiag(\0_{3\times 3}, I_{3(n-1)})]P^{-1}$ which implies that $A$ and $N$ have the same eigenvectors and therefore 
\begin{equation}\label{eq:hatA}
\hat A = NA = P[\blockdiag(\0_{3\times 3},\hat J_2)]P^{-1}.
\end{equation}
Similarly, it is easy to see that the Jordan normal form of $\tilde A$ is given by 
\begin{equation}\label{eq:tildeA}
\tilde A = \tilde P \hat J_1 \tilde P^{-1}.
\end{equation}}
\changed{
Therefore, under the conditions of Lemmas \ref{lem:multiplicity} and \ref{lem:jordan_chains} $\delta z_k$ and $\tilde z_k$ capture the behavior of a complementary set of eigenvalues of $A$. We are now ready to state our main convergence result.}


\begin{theorem}[Convergence]\label{th:convergence}
The following three statements are equivalent:
\begin{enumerate}
\item  The graph $G(V,E)$ is connected, $\kappa_1\neq\kappa_2$, $2>p>0$  and $\rho(\hat J_2)<1$ 
\item Condition \eqref{eq:synchronization_split} is satisfied
\item The algorithm \eqref{eq:system_z} achieves synchronization, i.e. \eqref{eq:synchronization} holds.
\end{enumerate}
Moreover, whenever the system synchronizes, we have
\begin{subequations}\label{eq:xs&rs}
\begin{eqnarray}\label{eq:xs}
x^* = \gamma\sum_{i=1}^n \xi_i\left(\frac{1}{r_i} x_{i}(t_0)+\tau\frac{\kp_2}{p^2}y_i(t_0)\right), \text { and }
\end{eqnarray}
\vspace{-.5cm}
\begin{eqnarray}\label{eq:rs}
r^* = \gamma\sum_{i=1}^{n} \xi_i(s_{i}(t_0) - \frac{\kappa_2}{p}y_{i}(t_0)).
\end{eqnarray}
\end{subequations}
\end{theorem}

The proof of Theorem \ref{th:convergence} can be found in Appendix \ref{app:th:convergence}. Theorem \ref{th:convergence} provides an analytical tool to understand the influence of the different nodes of the graph in the final offset $x^*$ and frequency $r^*$. For example, suppose that we know that node 1 has perfect knowledge of its own frequency ($r_1$) and the UTC time at $t=t_0$ ($x_1(t_0)=t_0$), and configure the network such that node 1 is the unique leader like the top node in Figures \ref{fig:topologies}a and \ref{fig:topologies}c. It is easy to show that $\xi_1=1$ and $\xi_i=0$ $\forall i\neq1$.  Then, using \eqref{eq:xs}-\eqref{eq:rs} and definition of $\gamma$ we can see that $\gamma=r_1$ and
\[
x^* = x_1(t_0) + r_1\tau \frac{\kp_2}{p^2}y_1(t_0)\text{ and }
r^* = r_1s_{1}(t_0) - \frac{r_1\kappa_2}{p}y_{1}(t_0).
\]
However, since node 1 knows $r_1$ and $t_0$, it can choose $x_1(t_0)=t_0$, $s_1(t_0)=\frac{1}{r_1}$ and $y_1(t_0)=0$. Thus, we obtain $x^*=t_0$ and $r^*=1$ which implies by \eqref{eq:synchronization} that every node in the network will end up with $x_i(t) = t $. In other words, Theorem \ref{th:convergence} allows us to understand how the information propagates and how we can guarantee that every server will converge to the desired time.
Notice that the initial condition used for server 1 is equivalent to assuming that server 1 is a reliable source of UTC like an atomic clock for instance.

\subsection{Necessary and sufficient conditions for synchronization}\label{ssec:param_vals}

\begin{figure}
\centering
\includegraphics[width=\columnwidth]{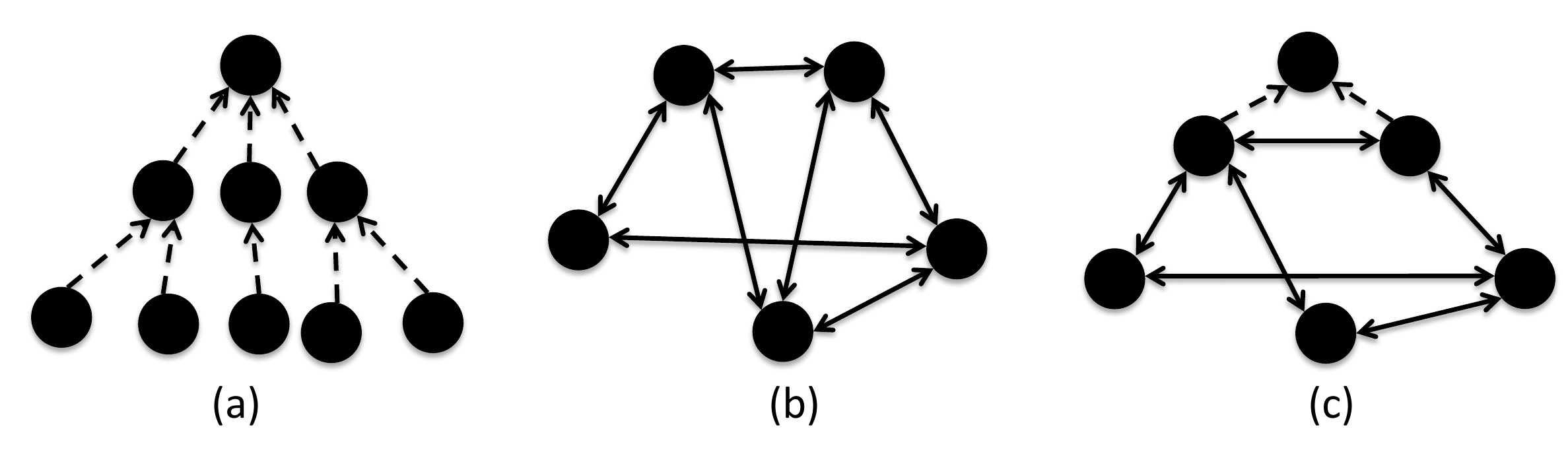}
\caption{Graphs with real eigenvalue Laplacians }\label{fig:topologies}
\end{figure}

We now provide necessary and sufficient conditions in terms of explicit parameter values ($\kp_1$, $\kp_2$ ,$\tau$ and $p$) for Theorem~\ref{th:convergence} to hold.
We will restrict our attention to graphs that have Laplacian matrices with real eigenvalues.
This includes for example trees (Figure \ref{fig:topologies}a), symmetric graphs with $\alpha_{ij}=\alpha_{ji}$ (Figure \ref{fig:topologies}b) and symmetric graphs with a leader (Figure \ref{fig:topologies}c).


The proof consists on studying the Schur stability of $g_l(\lambda)$ and has several steps. We first perform a change of variable that maps the unit circle onto the left half-plane. This transforms the problem of studying the Schur stability into a Hurwitz stability problem which is solved using  Hermite-Biehler Theorem which says:
Given the polynomial $P(s)=a_ns^n +...+a_0$, let $P^r(\omega)$ and $P^i(\omega)$ be the real and imaginary part of $P(j\omega)$, i.e. $P(j\omega) = P^r(\omega) +jP^i(\omega)$. Then $P(s)$ is a Hurwitz polynomial if and only if
\begin{enumerate}
\item $a_na_{n-1}>0$ and 
\item The zeros of $P^r(\omega)$ and $P^i(\omega)$ are all simple and real and interlace as $\omega$ runs from $-\infty$ to $+\infty$.
\end{enumerate}

\begin{theorem}[Parameter Values for Synchronization]\label{th:param_sync}
\changed{Consider a connected graph $G(V,E)$ with real eigenvalue Laplacian matrix $L$. Then, the system \eqref{eq:system_z} achieves synchronization 
if and only if}
\begin{enumerate}[(i)]
\item $|1-p| <1$ or equivalently $2>p>0$
\item $\frac{2\kp_1}{3p} >\kappa_1-\kappa_2 >0$ and (iii)
$\tau< \frac{p(\kappa_2 -p(\kappa_1-\kappa_2))}{\mu_{\max}(\kappa_1-p(\kappa_1-\kappa_2))^2}$
\end{enumerate}
where $\mu_{\max}$ is the largest eigenvalue of $LR$.
\end{theorem}

The proof of Theorem \ref{th:param_sync} can be found in Appendix \ref{app:th:param_sync}. Note that although $\mu_{\max}$ depends on $r_i$ which is in general unknown, it is easy to show that $\mu_l(LR)\leq\hat r_{\max}\mu_l(L)$  where $\hat r_{\max}$ is an upper bound of the maximum rate deviation $r_i$. Furthermore, using Greshgorin's circle theorem, it is easy to show that $\mu_{\max}(L)\leq2\alpha_{\max} :=2 \max_i \alpha_{ii}$.
Therefore, if we set
\begin{equation}\label{eq:tau_bound}
\tau <  \frac{p(\kp_2 - \delta\kp p)}{ 2\alpha_{\max}\hat r_{\max}(\kp_1-\delta\kp p)^2}
\end{equation}
convergence is guaranteed for {\bf every connected graph} with real  eigenvalues.

\section{Network Delays and Clock Wander}\label{sec:noise}

\changed{
In the previous section we showed the in the absence of network delays and clock wander, the system was able to achieve synchronization on a wide variety of communication topologies. In other words, we assumed the internal clock skew $r_i$ was fixed and that each computer could measure its offset with a neighbor $D_{ij}^x(t_k)$ without incurring in any error. We now study the behavior of our system when such assumptions are no longer true. We will model both, network delays and clock drifts using noise processes.
}

\changed{{\it Network Delays:}
Since our algorithm does not perform skew estimation, the network errors only affect the offset measurements $D^x_{ij}(t_k)$ in \eqref{eq:system}.} We use $g_{ij}^ww_{ij}(t_k)$ to denote the error incurred in estimating the offset between nodes $i$ and $j$ at time $t_k$, \changed{ i.e. we replace $D^x_{ij}(t_k)$ with $D^x_{ij}(t_k)+g_{ij}^ww_{ij}(t_k)$ in \eqref{eq:system}. This can be produced for instance by a congested connection between the two different nodes or due to path delay asymmetries~\cite{Veitch:2009td}.} We assume that $w_{ij}(t_k)$ has stationary mean $E\left[w_{ij}(t_k)\right]=\bar w_{ij}$ $\forall t_k$ and unit variance $E[(w_{ij}(t_k)-\bar w_{ij})^2]=1$ and use $g_{ij}^w$ to weigh the different connections. 

\changed{{\it Clock Wander:} We model the clock wander as a stochastic input to the clock skew adaptation. That is, instead of \eqref{eq:double_int_ib}, $s_i(t_k)$ changes according to 
\begin{equation}\label{eq:wander}
s_i(t_{k+1}) = s_i(t_k)  + u_i^s(t_k) + g^d_id_i(t_k)
\end{equation}
where $d_i(t_k)$ is a random variable with zero mean $E[d_i(t_k)]=0$ and unit variance $E[d_i(t_k)^2]=1$ and $g^d_i$ a positive scalar use to model clock heterogeneity. Equation \eqref{eq:wander} can be derived from a linear approximation of \eqref{eq:double_int} with a time varying internal skew $r_i(t_k)$ driven by an auto regressive process~\cite{Kim:2012db}, i.e.
 $ 
r_i(t_{k+1}) =  (1-q_i)r_i + q_ir_i(t_k) + g^d_id_i(t_k).
$ 
We omit the details of the derivation due to space constraints.}

This motivates the study of the stochastic process
\begin{align}\label{eq:system_noise}
z_{k+1}& = A z_k + Be_k,&
v_{k}& = Cz_k
\end{align}
where $e_k=[w_k^T \;d_k^T]^T$, $B=[ B_w\;B_d]$ with
{\small
\begin{equation}\label{eq:Bs}
B_w = \left[ \begin{array}{c}  \0_{n\times m} \\ -\kp_1 B_G^{-}\diag[\alpha_{ij}g_{ij}^w]\\-pB_G^{-}\diag[\alpha_{ij}g_{ij}^w]    \end{array}\right],\;
B_d = \left[ \begin{array}{c}  \0_{n\times n} \\ \diag[\beta_{i}g_{i}^d]\\  \0_{n\times n}     \end{array}\right],
\end{equation}
}$B_G^{-} = \min\{B_G,\0_{n\times m}\}$ and $B_G$ being the incidence matrix of $G(V,E)$ \footnote{ Notice that using this definition, we have $L=B_G^{-}\diag[\alpha_{ij}]B_G^T$.}  and  $w_k = [w_{ij}(t_k)]^T$.
The matrix $C$ maps the system state $z_k$ to the performance metric $v_k$ and will be specified in Section \ref{sec:experimental_results}.

\changed{One interesting difference between network delays ($w_{ij}(t_k)$) and wander ($d_i(t_k)$) is that in order to obtain good performance the algorithm should reject the noise from network delays $w_{ij}(t_k)$, but compensate the skew fast enough to follow $d_i(t_k)$.} 

In the remaining of this section, we first study the effect of biased network noise ($\bar w_{ij}\neq 0$) in the asymptotic frequency of the system and time offsets. In particular, we show that for arbitrarily distributed noise with stationary mean, the system's frequency tends to constantly drift unless there is a well defined leader in the topology. We then proceed to study how the parameters and network topology affect the systems performance, which is represented by the output signal $v_k$ of the stochastic process.





\subsection{Frequency Drift and Time Offsets}

\changed{Here, we concentrate on studying the evolution of the first moment of the stochastic process \eqref{eq:system_noise}. 
That is, we want to understand how $E[z_k]$ evolves as $k\rightarrow +\infty$. This is equivalent to study \eqref{eq:system_noise} in the case when the noise input $e_k$ is constant $e_k = \bar e=[\bar d^T \; \bar w^T]^T$, where $\bar d=\0_n$ by definition.}

\changed{Therefore, we will focus on understanding the effect of a constant input $\bar e$ on \eqref{eq:system_noise}. Again, we will use $\tilde z_k$ to understand how $\bar e$ affects the collective behavior of the clocks and $\delta z_k$ to understand how each individual clock deviates from the collective.
The next theorem summarizes the effect of non-zero mean error on the collective behavior of the system.}

\begin{theorem}[Frequency Drift]\label{th:frequency_drift}
\changed{In the presence of noise and under the condition of Theorem \ref{th:convergence} the collective frequency $\tilde s(t_k)$ will constantly drift away from its mean with probability one (in the set of possible $\bar w$), unless the graph $G(V,E)$ has a {\it unique leader}\footnote{A leader is a node to whom every other node can reach through a directed path.}. Whenever $G(V,E)$ does have a leader, the mean frequency $r^*$ is given by \eqref{eq:rs}.}
\end{theorem}


The proof of Theorem \ref{th:frequency_drift} can be found in Appendix \ref{app:th:frequency_drift_and_time_offsets}. 

\begin{remark}
\changed{Theorem \ref{th:frequency_drift} provides a precise characterization of how network delays and loops can produce instabilities. 
Similar results can be obtained for any protocol that controls clock speeds based on neighbors information. Therefore, this theorem shows that current industry practice is conservative and allows us to explore a wider set of topologies with timing loops (provided that such loops avoid the leader).}
\end{remark}


\changed{We now show how the deviations $\delta z_k$  are affected by $\bar e$.
\begin{theorem}[Time Offsets]\label{th:time_offsets}
Under the conditions of Theorem \ref{th:convergence} and graph $G(V,E)$ with unique leader, the deviations $\delta z_k$ converge to $\delta z^*$
given by
\[
\delta x^* = N_1 L^\dagger \delta w,\; \delta s^* = \0_n\text{ and } \delta y^* = \0_n.
\]
where $L^\dagger$ is the pseudo inverse of $L$ and $\delta w :=  -N_2B_G^{-}\diag[\alpha_{ij}g_{ij}^w]\bar w$.
\end{theorem}}

The proof can be found in Appendix \ref{app:th:frequency_drift_and_time_offsets}.

\subsection{$\mathcal H_2$ Performance Optimization}
We now proceed to study the effect of noisy measurements and wander on the output standard deviation of the system ($||v_k||_2$) when the input $e_k$ is white noise ($E[e_ke_l^T]=I_{m+n}\delta(l-k)$). In other words, we seek to find parameter values that minimize
\[
f(\kp_1,\kp_2,p,\alpha_{ij}) = ||v_k||_2 =\sqrt{ E\left[\lim_{N\rightarrow +\infty} \frac{1}{N}\sum_{k=0}^{N-1} v_k^Tv_k\right]}
\]

\changed{Since in practice we want to avoid any frequency drift introduced by the noise we will consider only topologies with a well defined leader.
Thus, all the randomness of the system is concentrated in $\delta z_k$ 
and we limit to study the stochastic process}
\begin{align*}
\delta z_{k+1} = \hat A\delta z + \hat Be_k,\quad v_{k} =  \hat C\delta z_k
\end{align*}
where $\hat A = NA$, $\hat B = NB$ and $\hat C = C$.



%
This optimization problem is standard in the control theory community and it can be shown to be equivalent to
\begin{subequations}\label{eq:H2}
\begin{align}
&\min_{X,\kp_1,\kp_2,p,\alpha_{ij}} \; f(\kp_1,\kp_2,p,\alpha_{ij}) \\
&\text{ subject to } \quad \rho(\hat A)\leq \rho^* \label{eq:rho_constrain}\\
&\quad X=\hat A^TX\hat A + \hat C^T\hat C
\end{align}
\end{subequations}
where $f(\kp_1,\kp_2,p,\alpha_{ij}):=\sqrt{\trace [X\hat B \hat B^T] }$, $\hat A$ is a function of $(\kp_1,\kp_2,p,\alpha_{ij})$ and $\rho^*<1$. The constrain \eqref{eq:rho_constrain} has been added in order to maintain the stability of $\hat A$.

\changed{While it is not in general easy to find the global minimum of \eqref{eq:H2}  there has been intensive research in designing optimization algorithms that find local minimums of the $\mathcal H_2$ norm of continuous time~\cite{arzelier_h2_2010} and discrete time~\cite{mostafa2008computational}  systems.
In this work  we solve \eqref{eq:H2} using a discrete-time version of the package Hifoo~\cite{arzelier_h2_2010,gumussoy_multiobjective_2009} known as Hifood~\cite{popov_fixed-structure_2010}. Several adaptations were needed to use Hifood to solve \eqref{eq:H2}. The details of these changes are documented in the Appendix \ref{app:Hifood}.}

\changed{The output of this optimization problem will be used in Experiment 6 of next section to demonstrate that the standard belief that ``clock precision degrades with the number of hops" is not necessarily true.}



\section{Experiments}\label{sec:experimental_results}

To test our solution and analysis, we implement  an asynchronous version of Algorithm 1 (Alg1) in C using the IBM CCT solution as our code base.
\changed{
Each server issues a thread to handle the connection with each neighbor. Every $\tau$ seconds (using OS time) each client takes offset measurements with its assigned neighbor and reports it to the main thread. Similarly, the main thread wakes up every $\tau$ seconds and gathers the offset information from all the connections and performs the update described in \eqref{eq:system}. We do not perform any explicit filtering of offset values, besides discarding spurious offsets larger than 500$ms$ in comparison with previous measurement.\footnote{An offset change of 500$ms$ within a $\tau$ of even 50 seconds implies a skew of  10,000ppm which is our maximum skew accepted.}
}

Our program reads the TSC counter directly using the {\verb rdtsc } assembly instruction to minimize reading latencies and maintains a virtual clock that can be directly updated. The list of neighbors is read from a configuration file and whenever there is no neighbor, the program follows the local Linux clock.  Finally, offset measurements are taken using an improved ping pong mechanism proposed in~\cite{froehlich_achieving_2008}.

\begin{figure}
\centering
\includegraphics[width=.8\columnwidth]{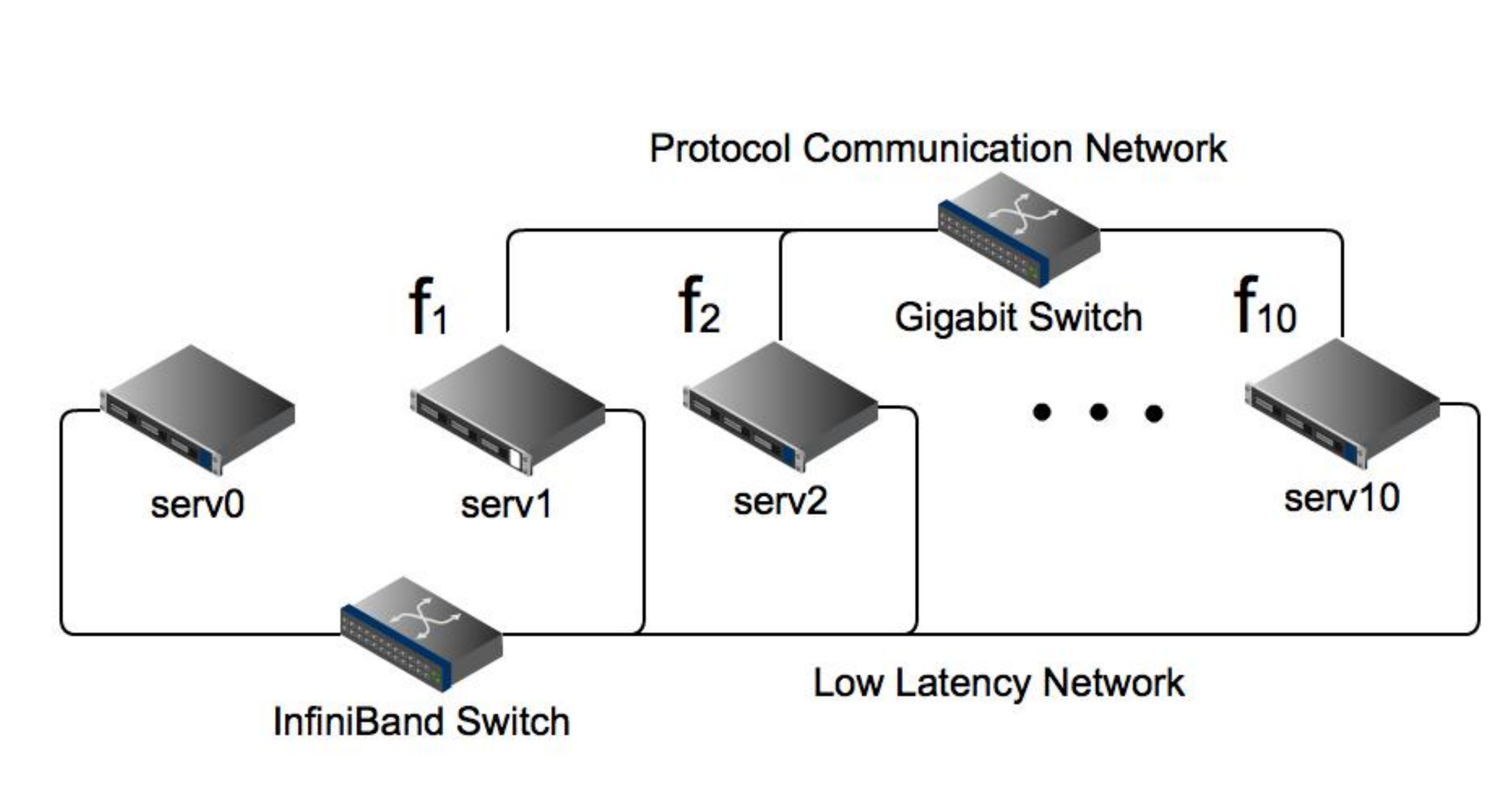}
\caption{Testbed of IBM BladeCenter blade servers} \label{fig:testbed}
\end{figure}

We run our skewless protocol in a cluster of IBM BladeCenter LS21 servers with two AMD Opteron processors of $2.40$GHz, and 16GB of memory. As shown in Figure \ref{fig:testbed}, the servers serv1-serv10 are used to run the protocol. The offset measurements are taken through a Gigabit Ethernet switch. 

\changed{Server serv0 is used as a common reference. It runs the same program that implements Alg1, but without skew adaptations, to measure the offset between itself and the other servers (serv1-serv10). These measurements are obtained through a 10Gbps Cisco 4x InfiniBand Switch to minimize network latencies.
Since the offset measurements performed by serv0 are done at different instances for different servers, we use linear interpolation to compensate this error. To compute the offset between two servers, say serv1 and serv2 ($x_1(t)-x_2(t)$), we subtract offset measurements obtained form serv0, i.e. $(x_1(t)-x_0(t))-(x_2(t)-x_0(t))$. Finally, we also eliminate spurious measurements that generate offsets bigger than 1ms as these are clearly due to network or os latencies.}

We use this testbed to validate the analysis in sections \ref{sec:analysis} and \ref{sec:noise}.
First, we illustrate the effect of different parameters and analyze
the effect of the network configuration on convergence (Experiment 1).
Then we present a series of configurations that demonstrate
how connectivity between clients is useful in reducing the
jitter of a noisy clock source (Experiment 2). We compare the
performance of our protocol with respect to NTP version 4 (Experiment 3) and 
\changed{IBM CCT (Experiment 4).  Finally, we verify the  presence frequency drift in the absence of a leader (Experiment 5), and study 
the interplay between network delays, wander and parameter values (Experiment 6).}

\changed{
The output performance signal $v_k$ will be the vector of offset differences between the leader $1$ and every other node $i$, i.e. $v_i(t_k)=x_i(t_k)-x_1(t_k)$ with $i\in\{2,...,n\}$. We will use a normalized version of it, referred here as {\it mean relative deviation} $\sqrt{S_n}$, as a performance metric. To make these comparison fair among different servers we correct our performance value by the empirical mean deviation to compensate biases due to path asymmetries.
In other words,
\begin{equation}
S_n= \frac{1}{{n-1}}{\sum_{i=2}^{n}\left<(x_i-x_1 -<x_i-x_1>)^2\right>}.
\end{equation}
where $<\cdot>$ amounts to the sample average.}
We will also use the $99\%$ Confidence Interval $CI_{99}$ and the maximum offset ($CI_{100}$) as metrics of accuracy. For example, if $CI_{99}=10\mu s$, then the $99\%$ of the offset samples will be within 10$\mu s$ of the leader.

Unless explicitly stated, the default parameter values are
\begin{align}\label{eq:values1}
p=0.99,\quad \kp_1 = 1.1, \quad  \kp_2 = 1.0 \text{ and }
\alpha_{ij}=\frac{c}{|\mathcal N_i|}.
\end{align}
The scalar $c$ is a commit or gain factor that will allow us to compensate the effect of $\tau$. Notice that by definition of $\alpha_{ij}$, $\alpha_{ii}=c$ for every node that is not the leader.


Moreover, these values immediately satisfy (i) and (ii) of Theorem \ref{th:param_sync} since $1-p=0.01$ and $\frac{2\kp_1}{3p} = 0.7407 > \kp_1-\kp_2 = 0.1$. The remaining condition can be satisfied by modifying $\tau$ or equivalently $c$. Here, we choose to fix $c=0.7$ which makes condition (iii)
\[
\tau< \frac{890.1}{\mu_{\max}}\text{ms}.
\]
For fixed polling interval $\tau$, the stability of the system depends on the value of $\mu_{\max}$, which is determined by the underlying network topology and the values of $\alpha_{ij}$.

\begin{figure}[htp]
\centering
\includegraphics[width=.55\columnwidth]{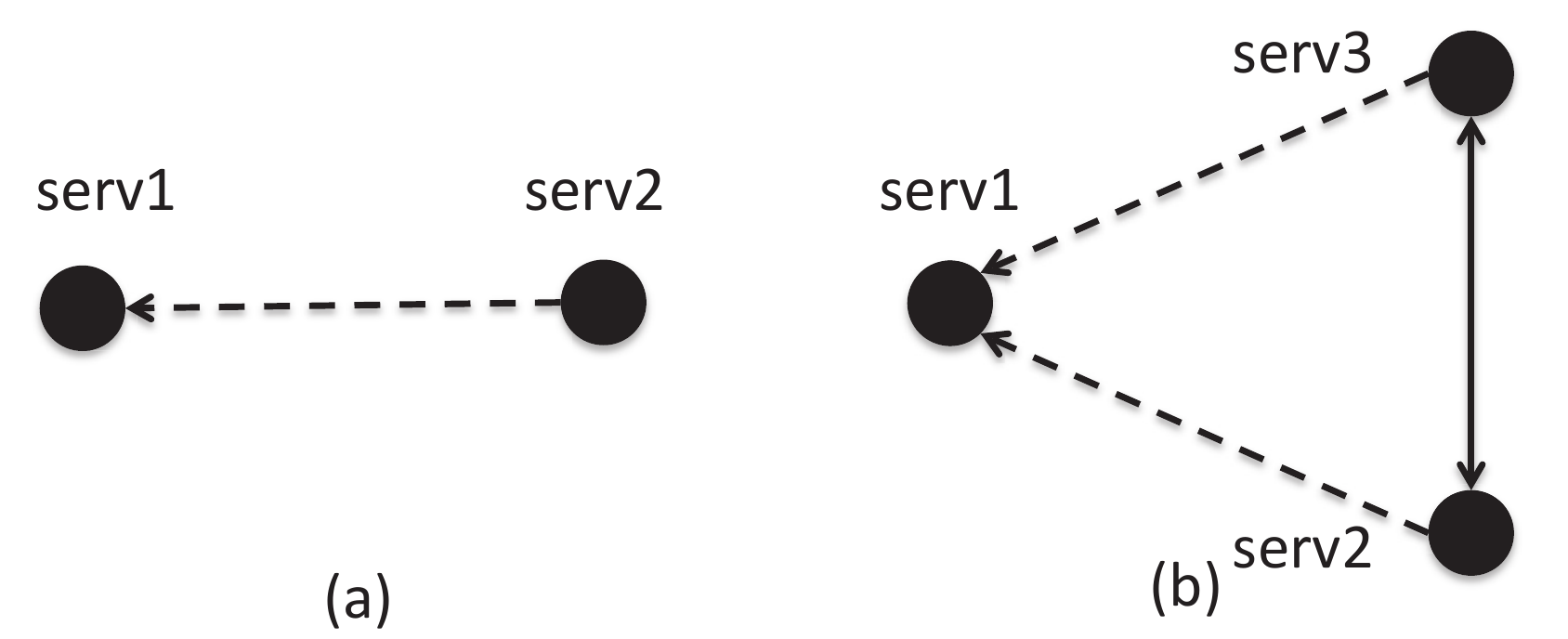}
\caption{Effect of topology on convergence: (a) Client-server configuration;  (b) Two clients connected to server and mutually connected.}\label{fig:23nodes}
\vspace{-.25cm}
\end{figure}

\begin{figure}[htp]
        \begin{subfigure}[b]{0.49\columnwidth}
               \centering
               \includegraphics[width=\columnwidth,]{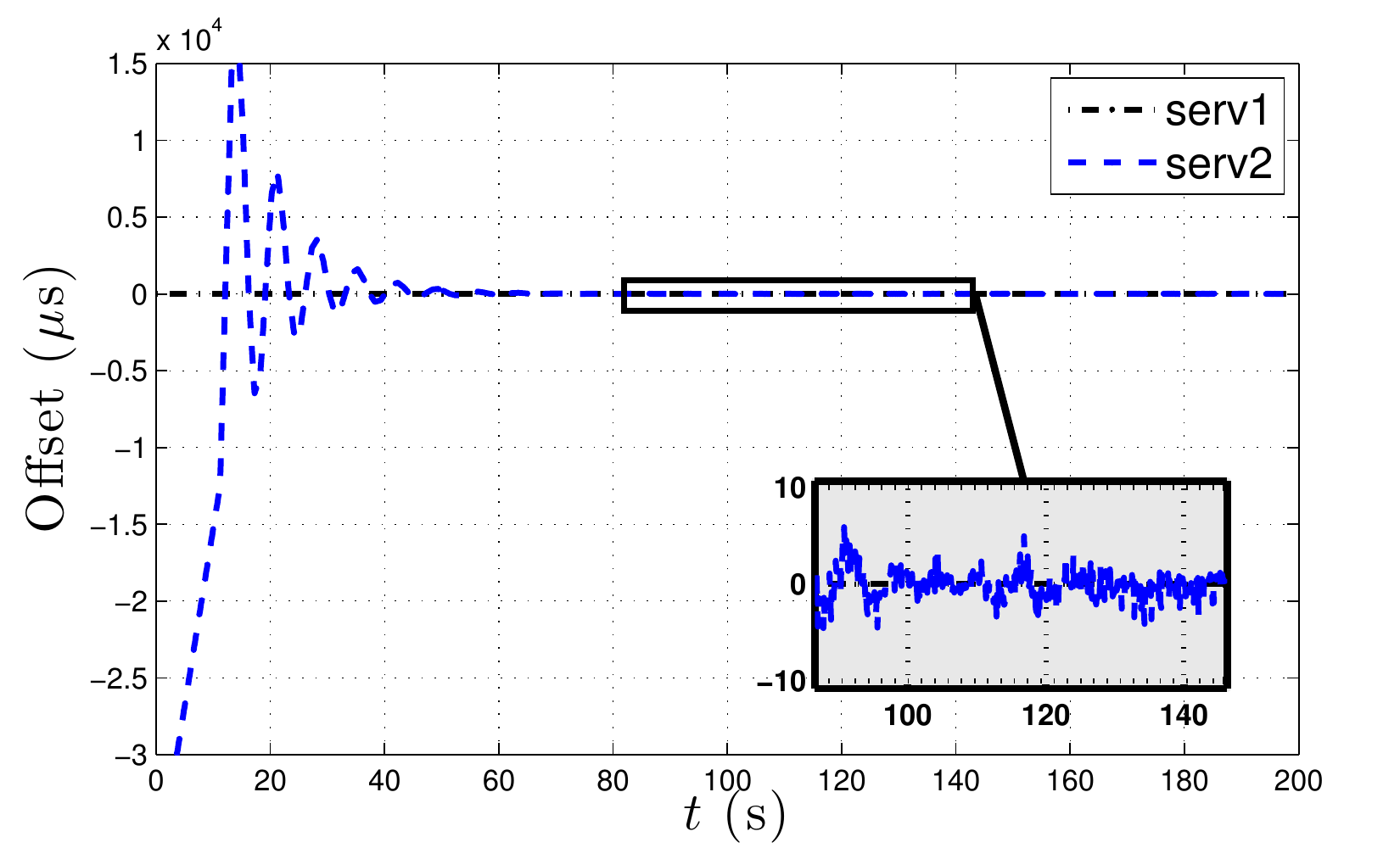}
               \caption{Client server configuration with $\tau = 1$s. The client converges and the algorithm is stable.}\label{fig:1000ms_21}
        \end{subfigure}%
        ~ 
        \begin{subfigure}[b]{0.49\columnwidth}
               \centering
               \includegraphics[width=\columnwidth]{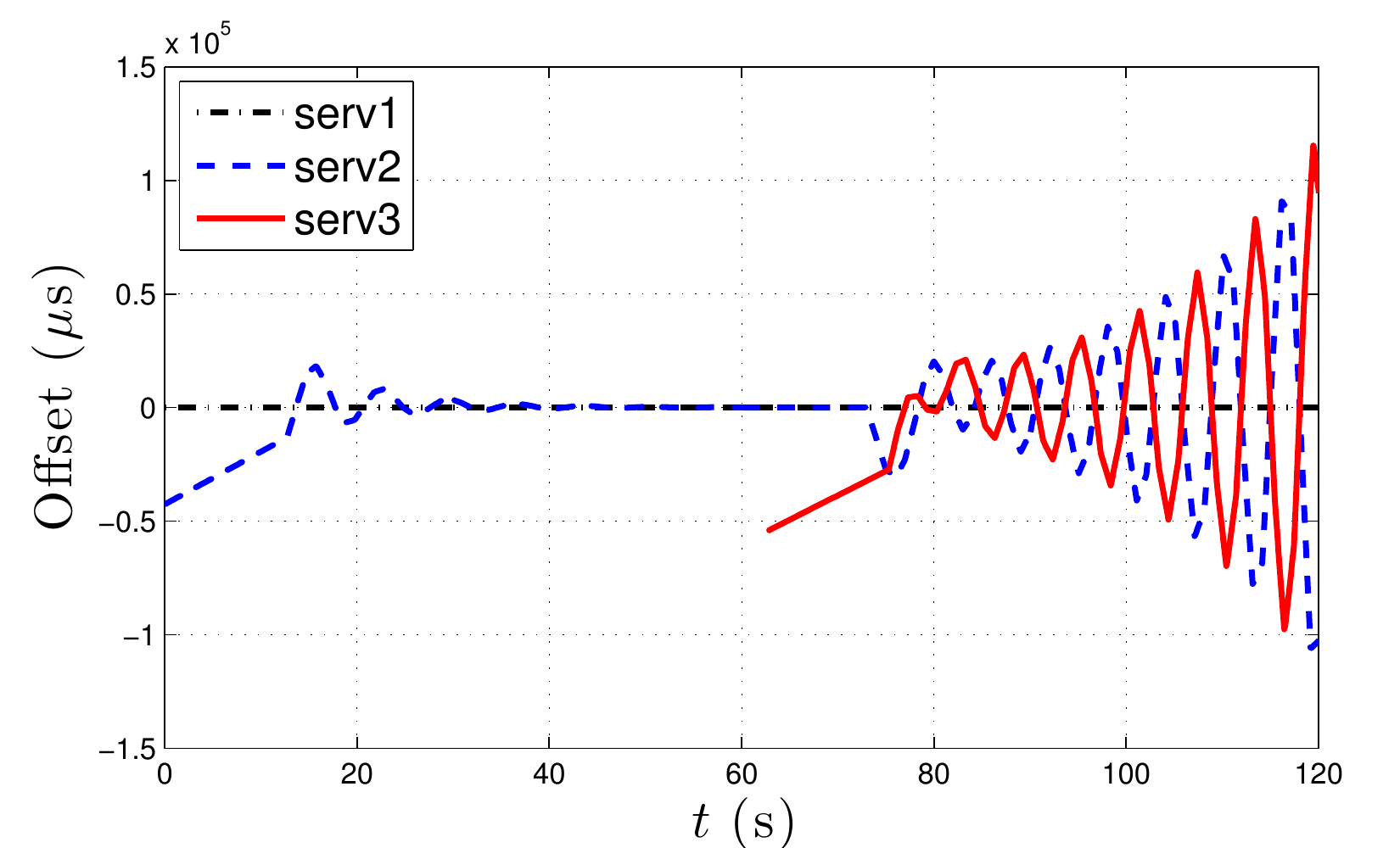}
               \caption{Two clients mutually connected with $\tau = 1$s. The algorithm becomes unstable.} \label{fig:1000ms_213}
        \end{subfigure}
        \caption{Loss of stability by change in the network topology}\label{fig:offsets}
\end{figure}

\vspace{.115cm}
\noindent{\bf Experiment 1 (Convergence):}
We first consider the client server configuration described in Figure \ref{fig:23nodes}a with a time step
$
\tau = 1\text{s}.
$ 
 In this configuration $\mu_{\max} \approx c=0.7$ and therefore condition (iii) becomes $\tau<1.2717$s. Figure \ref{fig:1000ms_21} shows the offset between serv1 (the leader) and serv2 (the client) in microseconds.
\changed{There we can see how serv2 gradually updates $s_2(t_k)$ until the offset becomes negligible.}

Figure \ref{fig:1000ms_21} tends to suggest that the set of parameters given by \eqref{eq:values1} and $\tau=1$s are suitable for deployment on the servers. This is in fact true provided that network is a directed tree as in Figure \ref{fig:topologies}a.
The  intuition behind this fact is that in a tree, each client connects only to one server. Thus, those connected to the leader will synchronize first and then subsequent layers will follow.

However, once loops appear in the network, there is no longer a clear dependency since two given nodes can mutually get information from each other. This type of dependency might make the algorithm unstable.
Figure \ref{fig:1000ms_213} shows an experiment with the same configuration as Figure \ref{fig:1000ms_21} in which serv2 synchronizes with serv1 until a third server (serv3) appears after $60$s. At that moment the system is reconfigured to have the topology of Figure \ref{fig:23nodes}b introducing a timing loop between serv2 and serv3. This timing loop makes the system unstable.

The instability arises since after serv3 starts, the new topology has  $\mu_{\max}\approx 1.5c=1.05$. Thus, the time step condition (iii) becomes
$
\tau< 847.8\text{ms}
$
which is no longer satisfied by $\tau=1$s.

\begin{figure}[htp]
\centering
\includegraphics[width=.85\columnwidth]{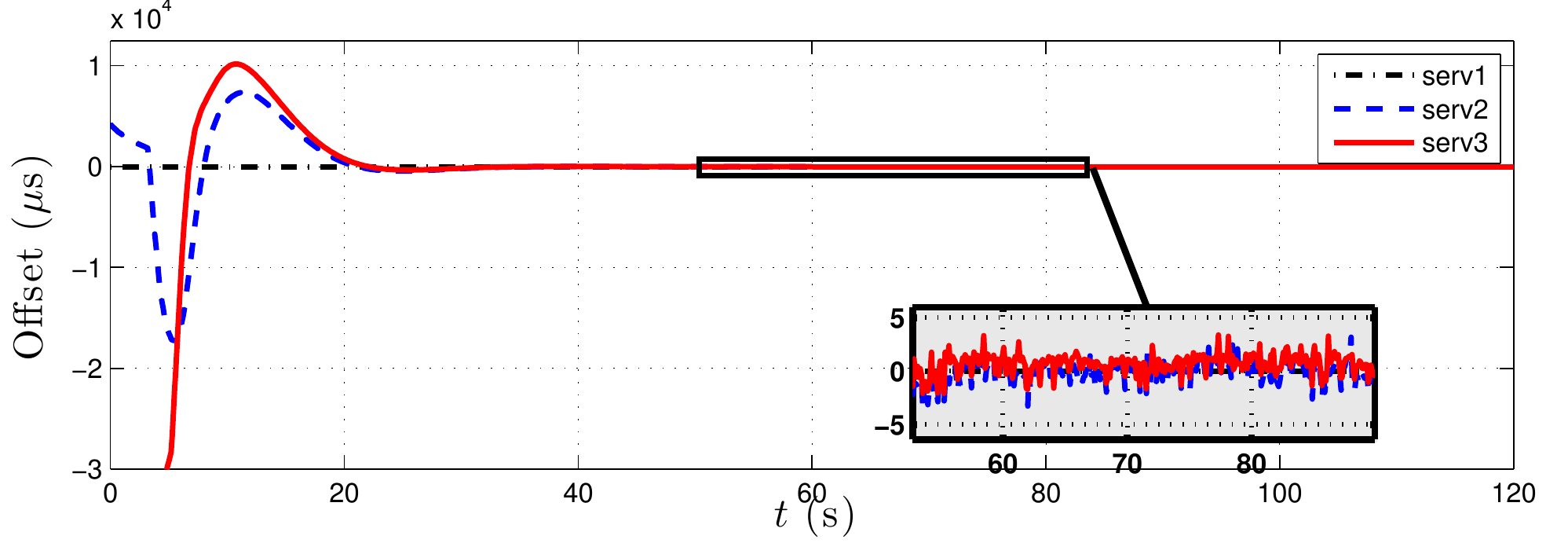}
\caption{\changed{Configuration of Figure \ref{fig:23nodes}b with $\tau = 500$ms. The algorithm becomes stable after reducing $\tau$ from $1$s to $500$ms.}}\label{fig:500ms_213}
\end{figure}

This may be solved for the new topology (Figure \ref{fig:23nodes}b) by using any $\tau$ smaller than $847.8$ms. However, if we want a set of parameters that is independent of the topology, we can use  \eqref{eq:tau_bound} and notice that $\alpha_{\max}=c$ and $\hat r_{\max}\approx 1$. We choose
\[
\tau = 500\text{ms} < \frac{890.2}{2\alpha_{\max}}\text{ms}=\frac{890.2}{2c}\text{ms}= 635.9\text{ms}.
\]
Figure \ref{fig:500ms_213} shows how  now serv2 and serv3 can synchronize with serv1 after introducing this change.

\noindent
{\bf Experiment 2 (Timing Loops Effect):}
 We now show how timing loops can be used to collectively outperform individual clients when the time source is noisy.

\begin{figure}[htp]
\centering
\includegraphics[width=.8\columnwidth]{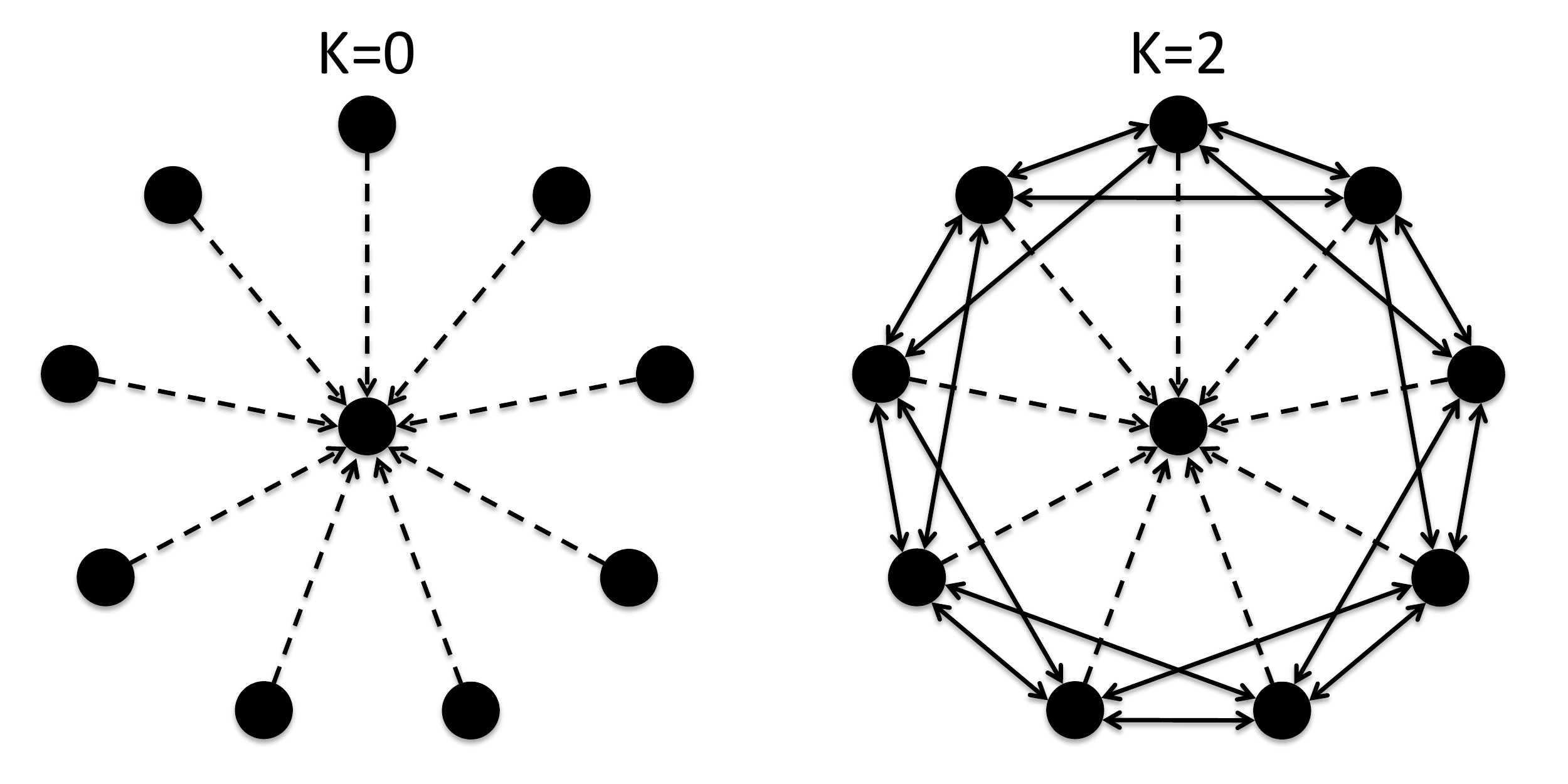}
\caption{Leader topologies with $2K$ neighbors connection. Connections to the leader (serv1) are unidirectional while the connections among clients (serv2 through serv10) are bidirectional}\label{fig:wheel}
\end{figure}

We run Alg1 on 10 servers (serv1 through serv10). The connection setup is described in Figure \ref{fig:wheel}. Every node is directly connected unidirectionally to the leader (serv1) and bidirectionally to $2K$ additional neighbors. 

When $K=0$ then the network reduces to a star topology and when $K=4$ the servers serv2 through serv10 form a complete graph.
\changed{The dashed arrows in Figure \ref{fig:wheel} show the connections where jitter was introduced. To emulate a link with jitter we added random noise $\eta$ with values taken uniformly from $\{0,1,...,\text{Jitter}_{\max}\}$ on both directions of the communication,}
\begin{equation}
\eta \in \{0,1,...,\text{Jitter}_{\max}\}\text{ms}.\label{eq:jitter}
\end{equation}

Notice that the arrow only shows a dependency relationship, the ping pong mechanism sends packets in both direction of the physical communication. We used a value of Jitter$_{\max}=10$ms.
Since the error was introduced in both directions of the ping pong, this is equivalent to a standard deviation of $6.05$ms.

 \begin{figure}[htp]
        \begin{subfigure}[b]{0.49\columnwidth}
                \centering
                \includegraphics[width=\columnwidth]{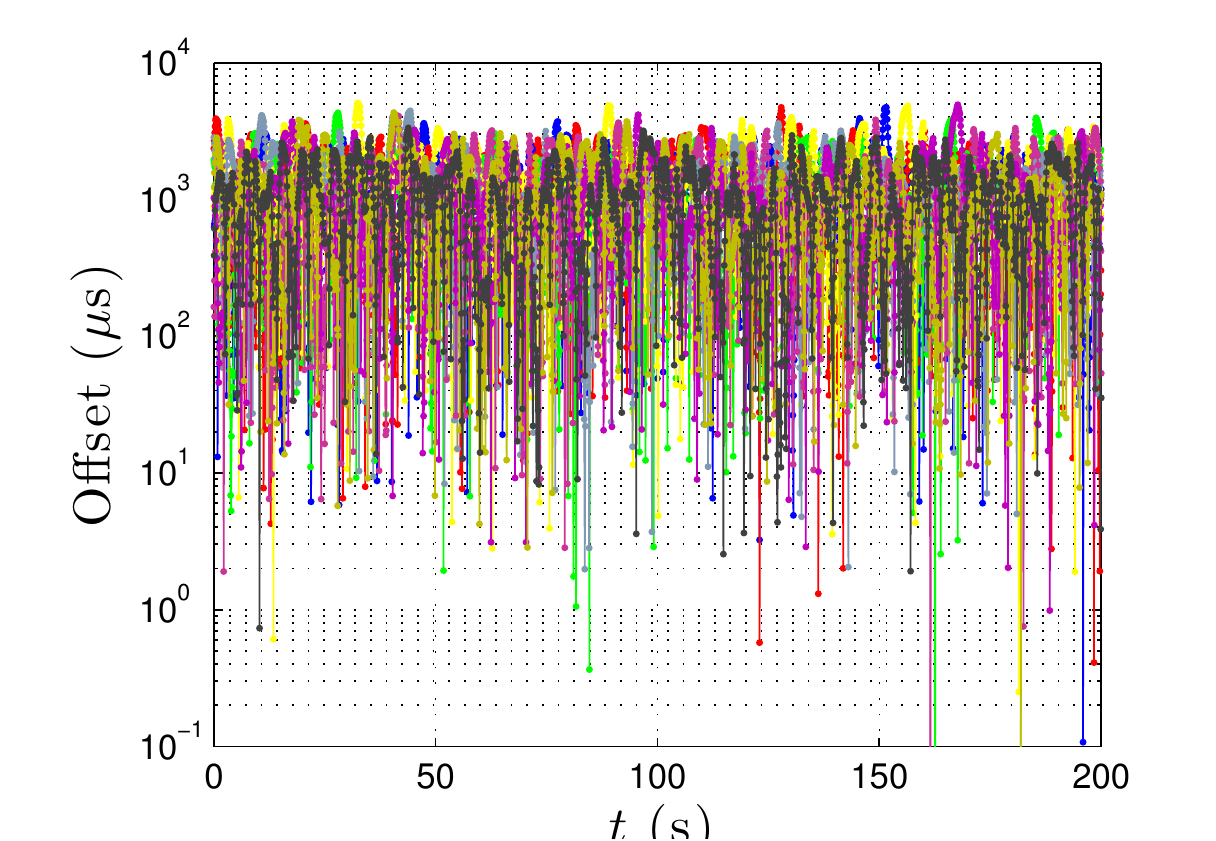}
                \caption{Star topology ($K=0$)}
                \label{fig:offset_star}
        \end{subfigure}%
        ~ 
        \begin{subfigure}[b]{0.49\columnwidth}
                \centering
                \includegraphics[width=\columnwidth]{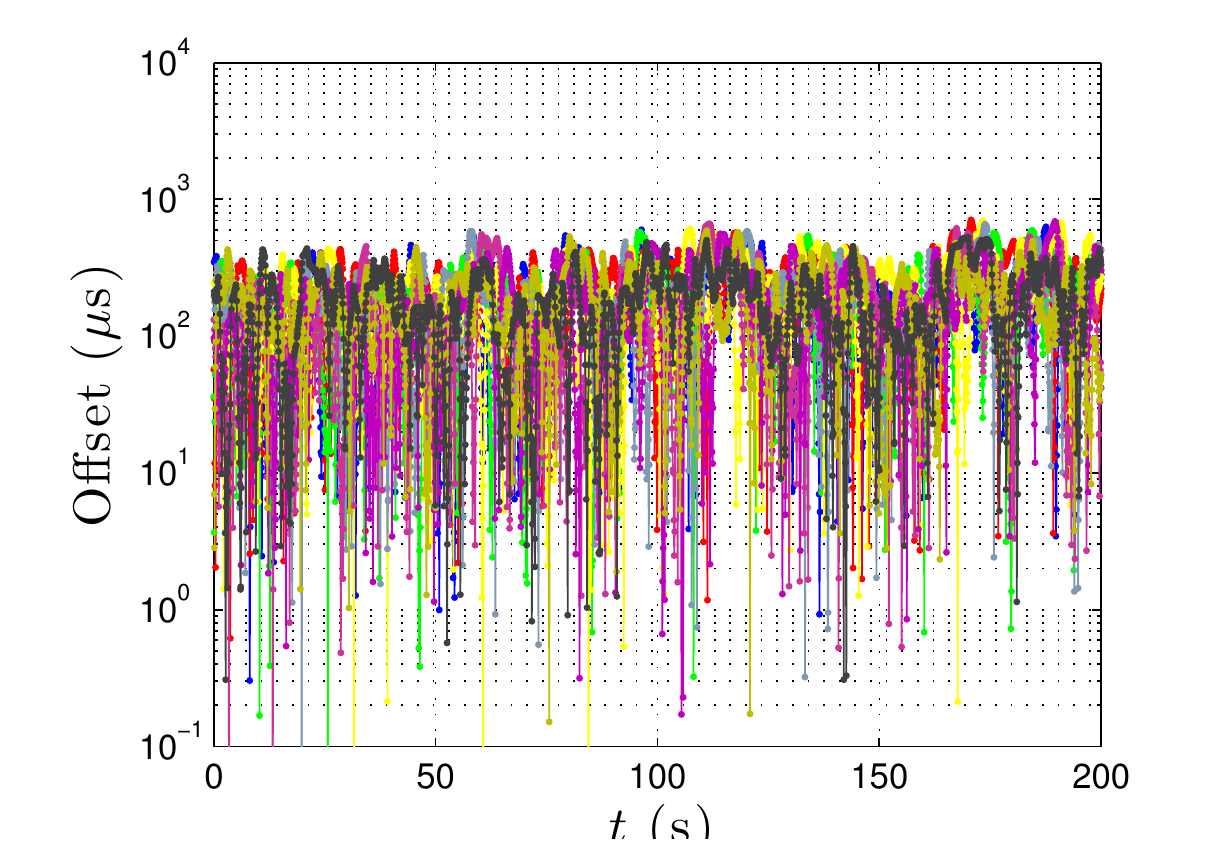}
                \caption{Complete subgraph ($K=4$)}
                \label{fig:offset_K}
        \end{subfigure}
        \vspace{-.25cm}
        \caption{Offset of the nine servers connected to a noisy clock source}\label{fig:offsets}
\end{figure}

Figure \ref{fig:offsets} illustrates the relative offset between the two extreme cases; The star topology ($K=0$) is shown in Figure \ref{fig:offset_star}, and the complete subgraph ($K=4$) is shown in Figure \ref{fig:offset_K}.

The worst case offset for $K=0$ is $CI_{100}=5.1$ms which is on the order of the standard deviation of the jitter. However, when $K=4$ we obtain a worst case offset of $CI_{100}=690.8\mu$s, an {\bf order of magnitude improvement.}

\begin{figure}[htp]
\centering
\includegraphics[width=.85\columnwidth]{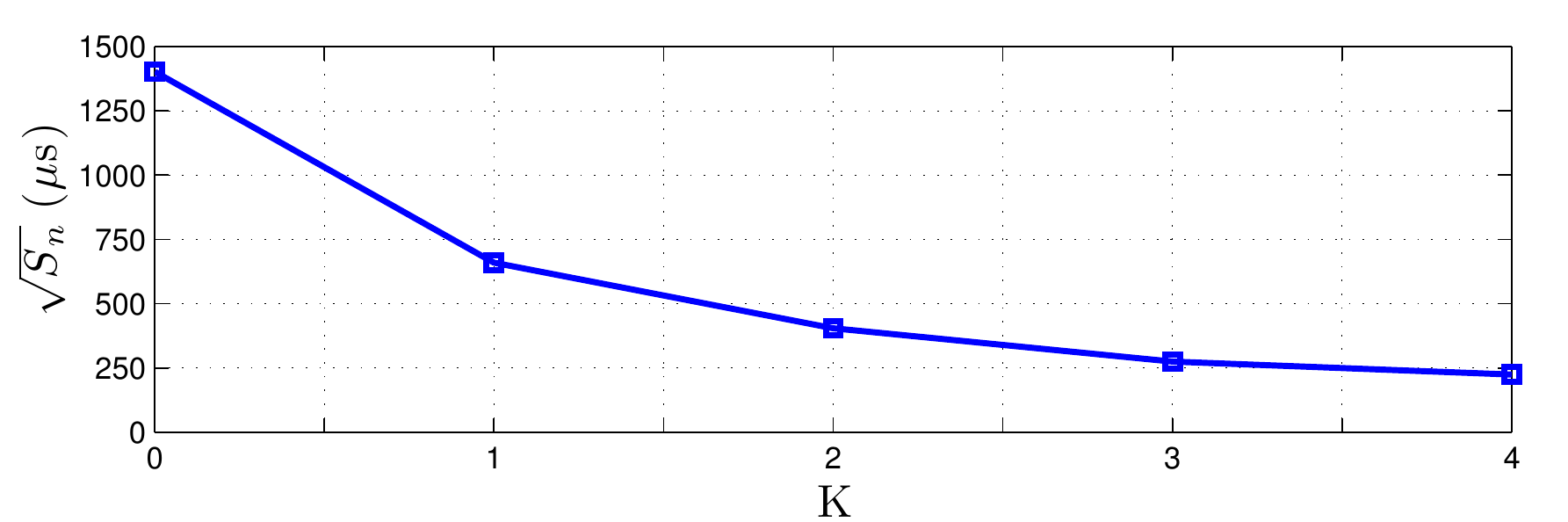}
\caption{Effect of the client's communication topology on the mean relative deviation. As the connectivity increases ($K$ increases) the mean relative deviation is reduced by factor of  $6.26$, i.e. a noise reduction of approx. 8dB.}\label{fig:error_wheel}
\end{figure}

The change on the mean relative deviation $\sqrt{S_n}$ as the connectivity among clients increases is studied in Figure \ref{fig:error_wheel}.
The results presented show that even without any offset filtering mechanism the network itself is able to perform a distributed filtering that achieves an improvement of up to a factor of $6.26$ in $\sqrt{S_n}$, or equivalently a noise reduction of almost 8dB.

\vspace{.115cm}
\noindent{\bf Experiment 3 (Comparison with NTPv4):}
We now perform a thorough comparison between our protocol (Alg1) and NTPv4.
We used a one hop configuration using serv1 as leader running an NTPv4 server and Alg1,  and serv9 and serv10 as clients, connected only to serv1, running NTPv4 and Alg1 respectively.

In order to make a fair comparison, we need both algorithms to use the same polling interval. Thus, we fix $\tau=16$sec. This can be done for NTP by setting the parameters  {\verb minpoll } and {\verb maxpoll } to $4$ ($2^4=16$secs).
\changed{The remainder parameter values for Alg1 were obtained with our optimization framework (using $g_{ij}^w=100$ and $g_{ij}^d=1e-3$) and are given by}
\begin{align}\label{eq:values3}
p=1.98,\quad \kp_1 = 1.388 \text{ and }  \kp_2 = 1.374.
\end{align}

\begin{figure}[htp]
       \begin{subfigure}[b]{0.49\columnwidth}
               \centering
               \includegraphics[height=.6\columnwidth,width=\columnwidth]{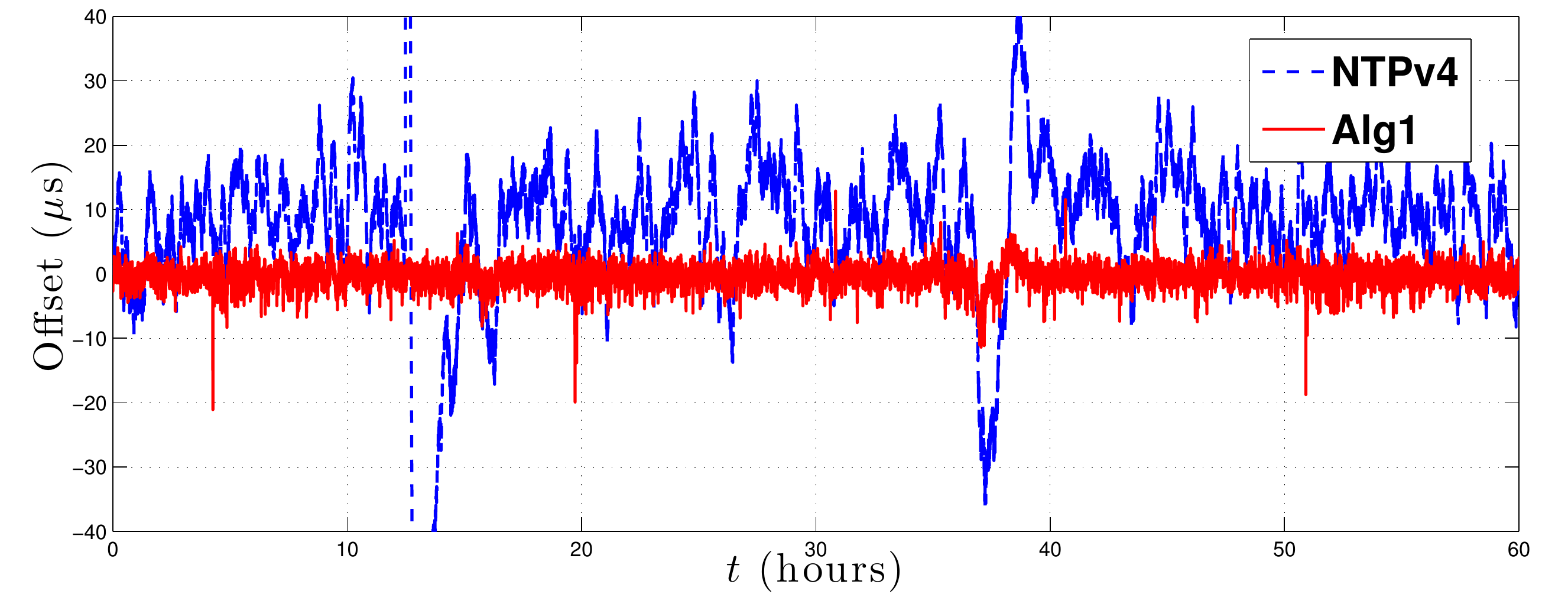}
               \caption{Offset values of NTPv4 and Alg1 for a period of 60 hours.\\}\label{fig:offset_ntpvsalg1}
        \end{subfigure}
        \begin{subfigure}[b]{0.49\columnwidth}
               \centering
			  \includegraphics[height=.6\columnwidth,width=\columnwidth]{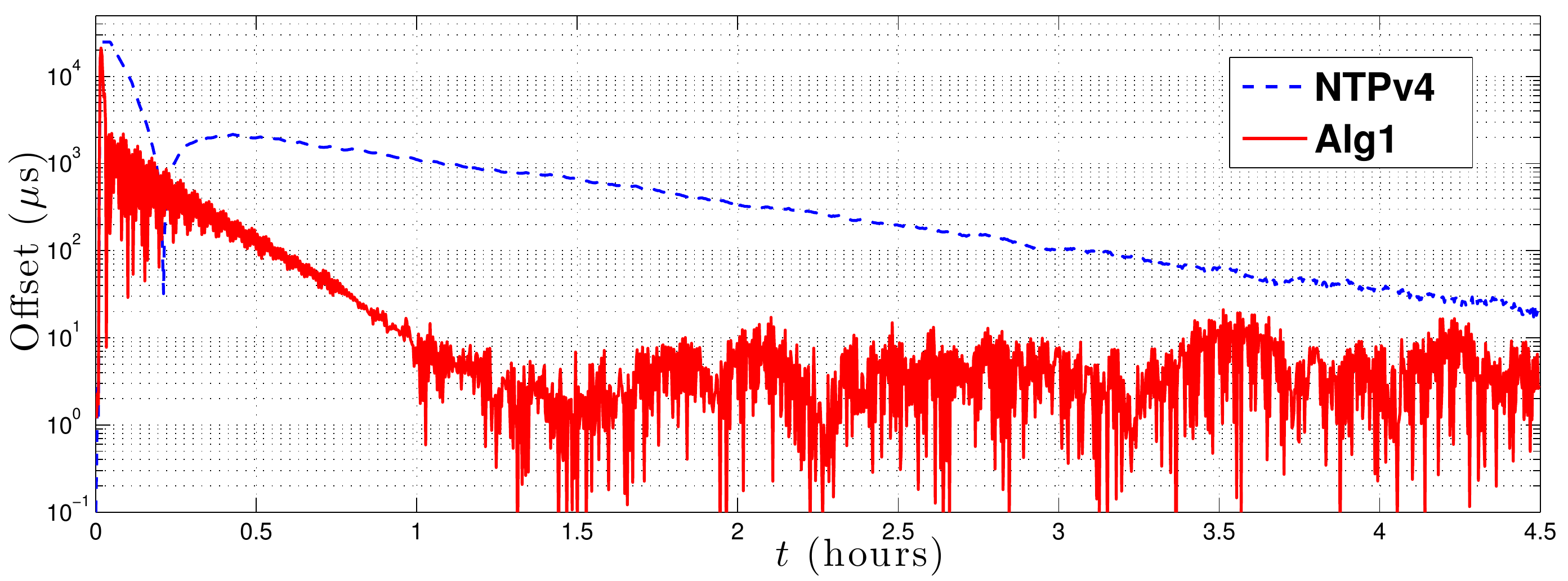}
			  \caption{Offset values of NTPv4 and Alg1 after a 25ms offset introduced in serv1. }\label{fig:step_ntpvsalg1}
        \end{subfigure}%
       \vspace{-.25cm}
        \caption{Performance evaluation between our solution (Alg1) and NTPv4}\label{fig:comparison_ntp}
\end{figure}

\changed{Figure \ref{fig:offset_ntpvsalg1} shows the time differences between the clients running NTPv4 and Alg1 (serv9 and serv10) , and the leader (serv1) over a period of 60 hours.  It can be seen that Alg1 is able to track serv1's clock keeping an offset smaller than 5$\mu$s for most of the time while NTPv4 incurs in larger offsets during the same period of time. This difference is produced by the fact that Alg1 is able to react more rapidly to frequency changes while NTPv4 incurs in more offset corrections that generate larger jitter.}

\changed{The mean offsets of Alg1 and NTPv4 are  $-0.48\mu$s and $9.00\mu$s. This difference in mean is mainly due to an asymmetric path on the measurements from serv9 to serv1. After compensating this bias, Alg1 achieves a performance of $\sqrt{S_n}=1.3\mu$s, $CI_{99}=4.9\mu$s and a maximum offset of $CI_{100}=20.6\mu$s, while NTPv4 obtains
$\sqrt{S_n}=6.4\mu$s, $CI_{99}=74.5\mu$s and a maximum offset of $CI_{100}=1.4$ms. 
Thus, not only Alg1 achieves a reduction of $\sqrt{S_n} $ by a factor of $5.0$ ($-7$dB) with respect to NTPv4, but it also obtains smaller confidence intervals and maximum offset values. 
}

\begin{figure}[htp]
               \centering
               \includegraphics[width=\columnwidth]{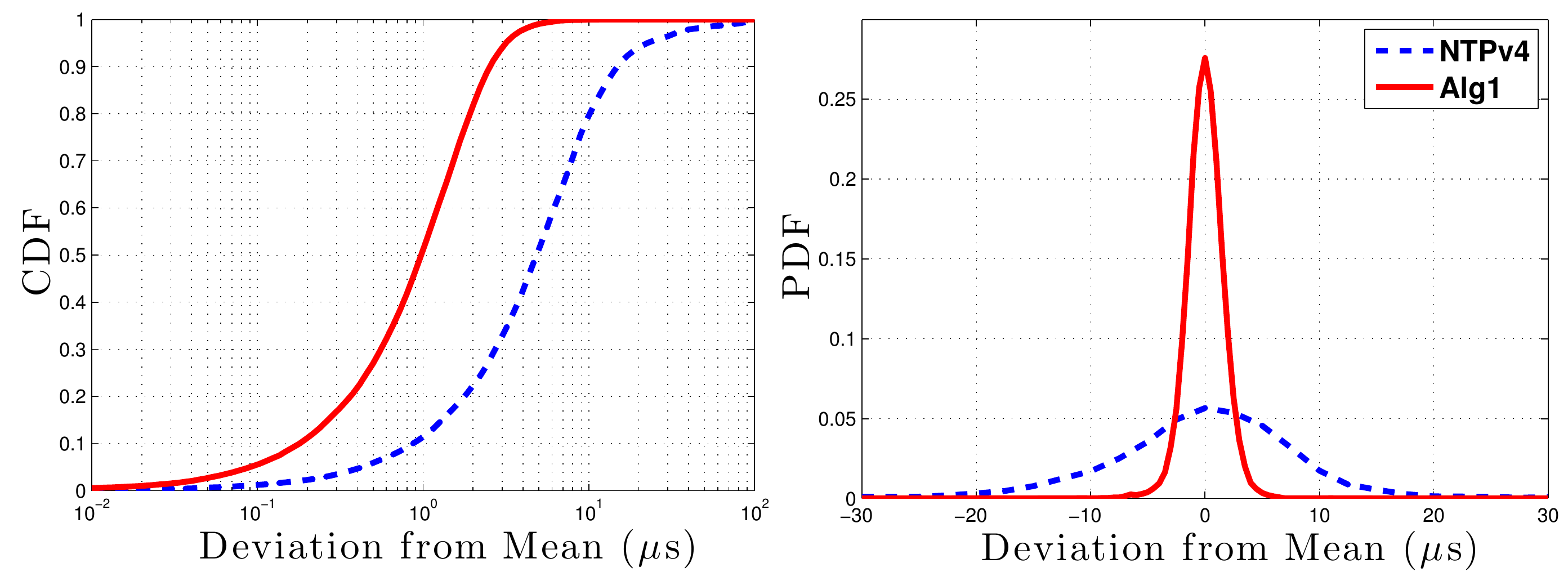}
        \caption{Empirical Cumulative Distribution Function (CDF) and Probability Density Function (PDF) of  Alg1 and NTPv4}\label{fig:cdf_pdf_ntpvsalg1}
\end{figure}

\changed{A more detailed and comprehensive analysis is presented in Figure \ref{fig:cdf_pdf_ntpvsalg1} where we plot the Cumulative Distribution Function (CDF) and Probability Density Function of the samples. 
The improvement of Alg1 with respect to NTPv4 is again clearly seen here.
}


%

Finally, we investigate the speed of convergence. Starting from both clients synchronized to server serv1, we introduce a 25ms offset. Figure \ref{fig:step_ntpvsalg1} shows how Alg1 is able to converge to a $20\mu s$ range within one hour while NTPv4 needs $4.5$hours to achieve the same synchronization precision.

\vspace{.115cm}
\noindent{\bf Experiment 4 (Comparison with IBM CCT):}
We now proceed to compare the performance of Alg1 with respect to IBM CCT. Notice that unlike IBM CCT, our solution does not perform any previous filtering of the offset samples, the filtering is performed instead by calibrating the parameters. Here we use  $c=0.70$, $\tau=250ms$, $\kp_1=0.1385$, $\kp_2=0.1363$ and $p=0.62$.

\begin{figure}[htp]
       \begin{subfigure}[b]{0.49\columnwidth}
               \centering
               \includegraphics[width=\columnwidth]{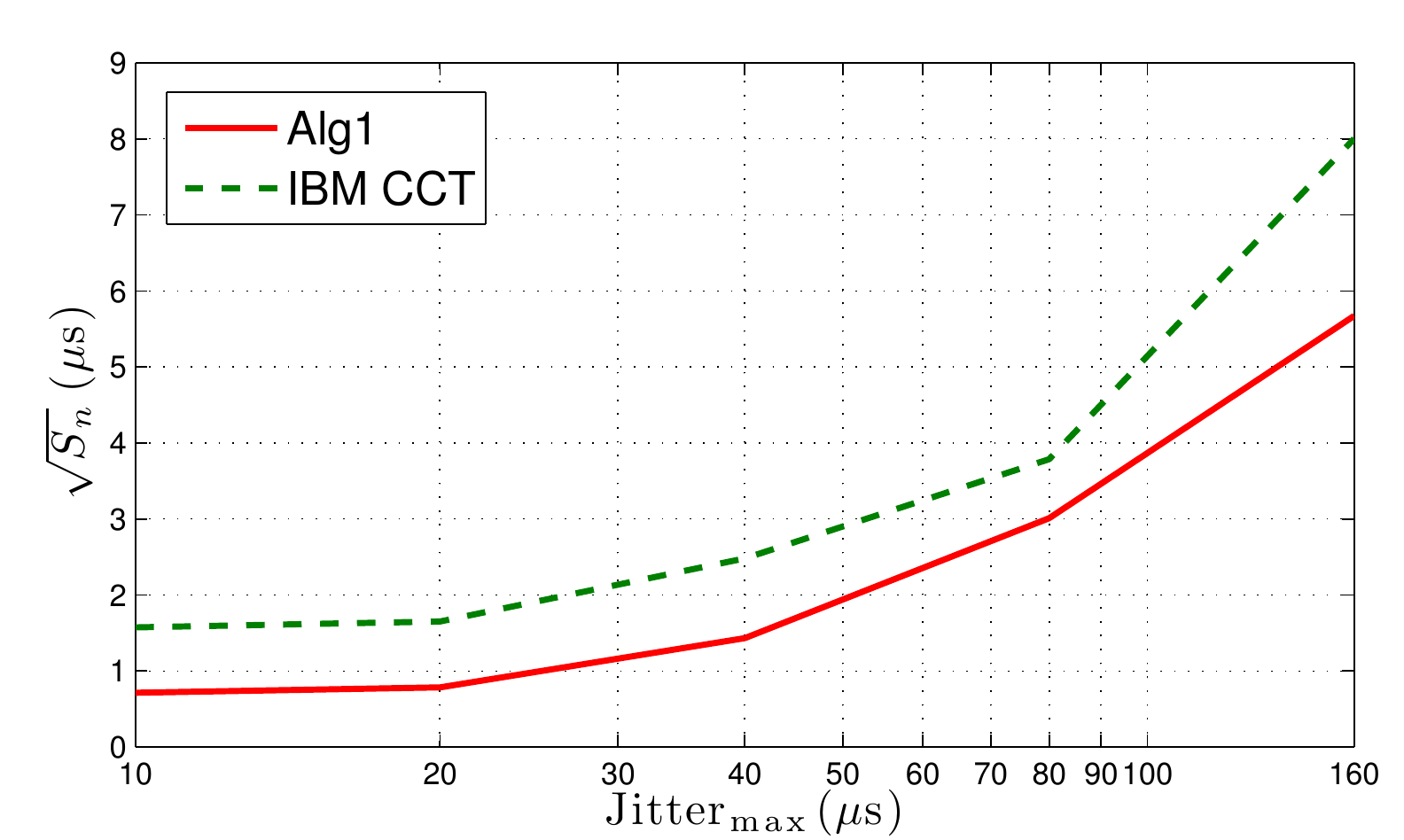}
               \caption{Mean relative deviation $\sqrt{S_n}$}\label{fig:comparison1}
        \end{subfigure}
        \begin{subfigure}[b]{0.49\columnwidth}
               \centering
               \includegraphics[width=\columnwidth]{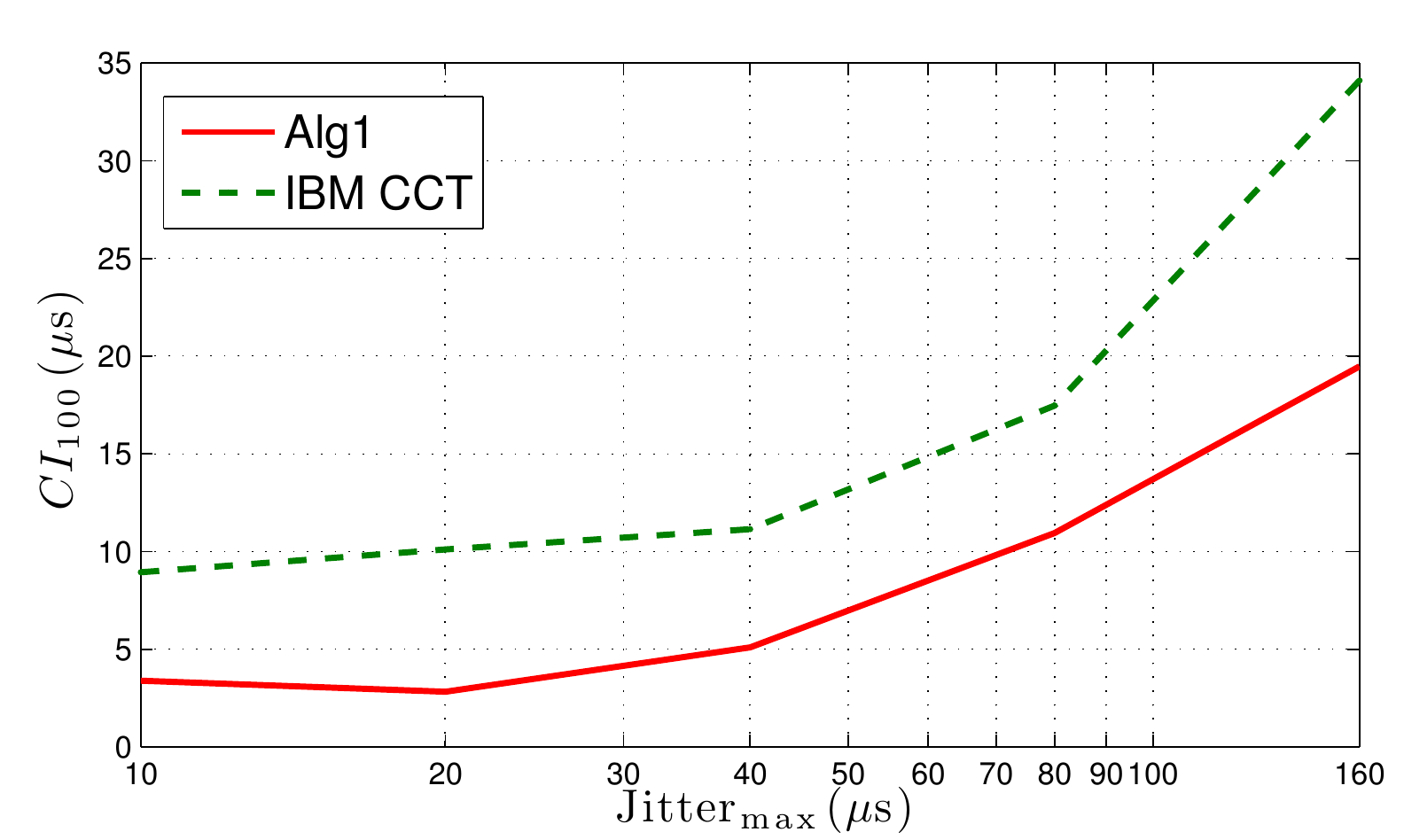}
               \caption{Maximum offset}\label{fig:comparison2}
        \end{subfigure}%
        \vspace{-.25cm}
        \caption{Performance evaluation between our solution (Alg1) and IBM CCT}\label{fig:comparison}
\end{figure}

In Figure \ref{fig:comparison1} we present the mean relative deviation $\sqrt{S_n}$ for two clients connected directly to the leader as the jitter is increased from Jitter$_{\max}=10\mu$s  to Jitter$_{\max}=160\mu$s, doubling Jitter$_{\max}$ each time, with a granularity in the random generator of $1\mu$s. The worst case offset is shown in Figure \ref{fig:comparison2}. Each data point is computed using a sample run of 250 seconds.

Our algorithm consistently outperforms IBM CCT in terms of both $\sqrt{S_n}$ and worst case offset.
The performance improvement is due to two reasons. Firstly, the noise filter used by the IBM CCT algorithm is tailored for noise distributions that are mostly concentrated close to zero with sporadic large errors. However, it does not work properly in cases where the distribution is more homogeneous as in this case.
Secondly, by choosing $\delta\kp = \kp_1-\kp_2=0.002\ll1$
the protocol becomes very robust to offset errors.

\vspace{.115cm}
\noindent{\bf Experiment 5 (Frequency drift without leader):}
We now proceed to experimentally verify that without leader, the system tends to constantly drift the frequency. Our analysis predicts that even the minor bias in the offset measurements will produce this effect. To verify this phenomenon, we use the network topology in Figure \ref{fig:23nodes}b with $\tau=0.5$s and wait for the system to converge.
\begin{figure}[htp]
\centering
\includegraphics[width=.9\columnwidth]{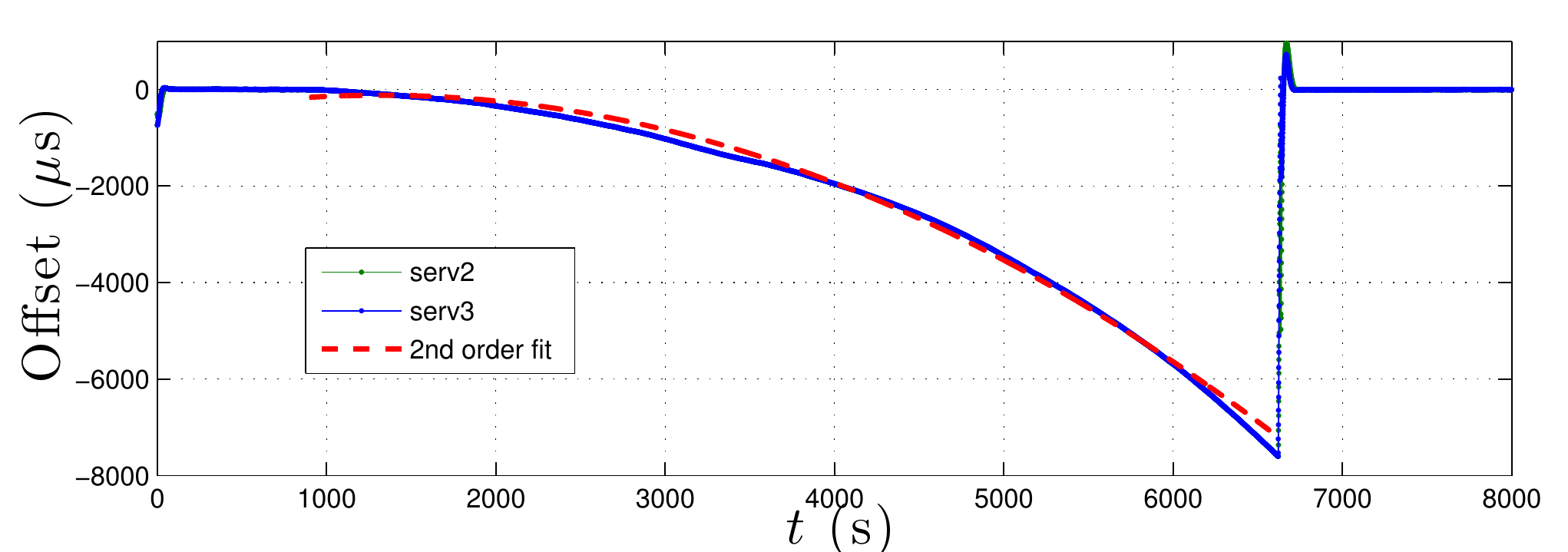}
\caption{Frequency drift}\label{fig:freq_drift}
\end{figure}

After 1000s the timing process of serv1 is turned off. Figure \ref{fig:freq_drift} shows how the offsets of serv2 and serv3 start to grow in a parabolic  trajectory characteristic  of a constant acceleration, i.e. constant drift. After 6600s serv1 is restarted and the system quickly recovers synchronization.
A second order fit of the faulty trajectory was perform obtaining a drift of approximately $-250\text{ ns}/\text{s}^2$. While this is not quite significant in the first few minutes, it becomes significant as time goes on.

\begin{figure}[htp]
\centering
\includegraphics[width=1.0\columnwidth]{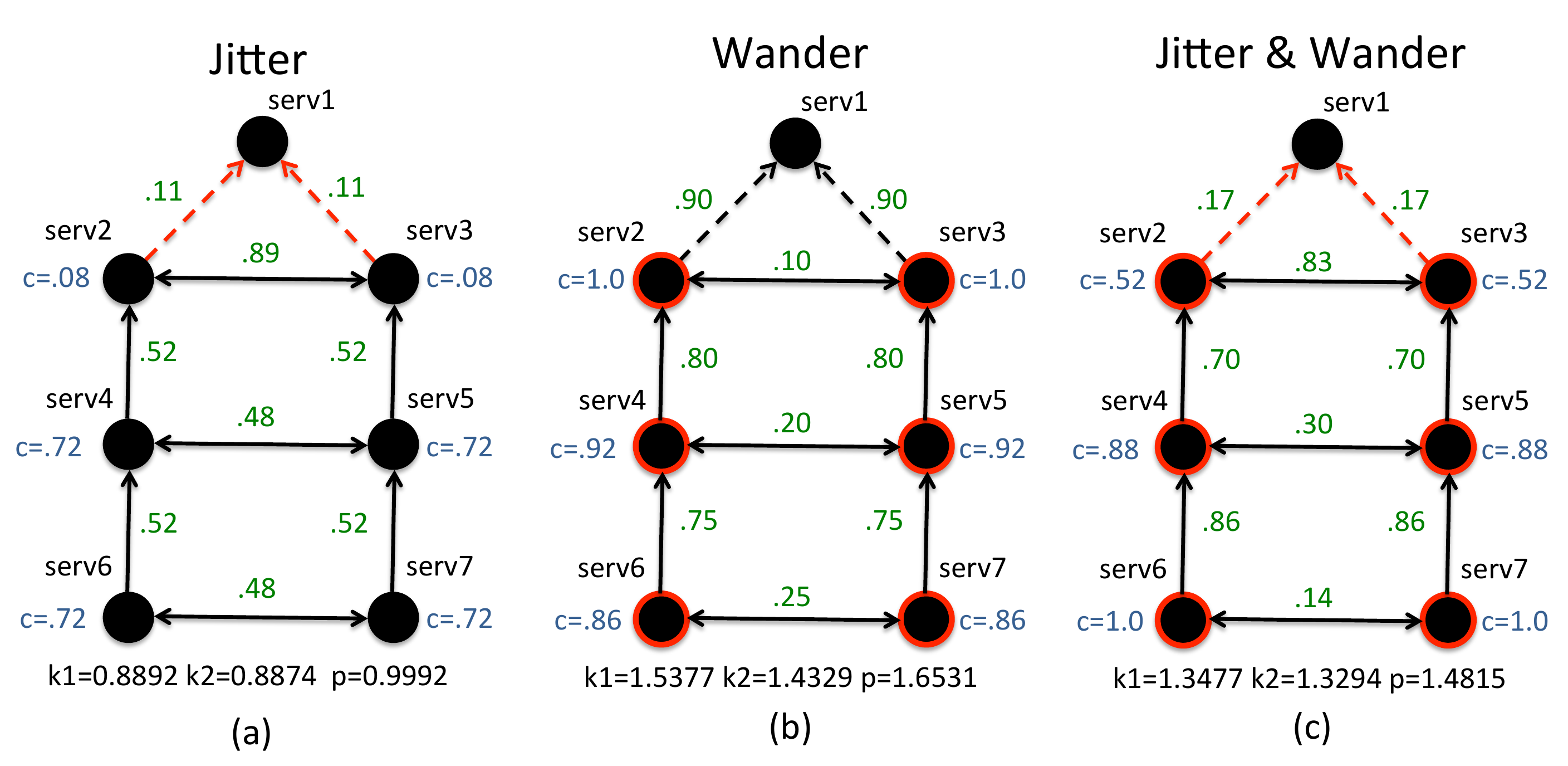}
\caption{Network scenarios and optimal parameters}\label{fig:exp_h2}
\end{figure}

\vspace{.115cm}
\noindent{\bf Experiment 6 (Jitter and Wander Tradeoff):}
Finally, we use the proposed $\mathcal H_2$ optimization scheme to show how the optimal parameter values depend on the different noise conditions within the network described in Figure \ref{fig:exp_h2}. We consider three different noise scenarios in which we either add jitter between server serv1 and servers serv2 and serv3, and/or add wander on severs serv2-serv7. In all the cases we used $\tau = 0.5$s and make offset measurements through the InfiniBand switch to minimize the any additional source of noise.

The jitter is generated  by adding in both directions of the physical communication a random value $\eta$ similarly to Experiment 2 (c.f. \eqref{eq:jitter}), but with a Jitter$_{\max}=100\mu$s. This generates an aggregate offset measurement noise of zero mean and standard deviation of $40.8\mu$s. On the other hand, the wander is generated by adding gaussian noise with zero mean and standard deviation of $0.2$ppm in the $s_i(t_k)$ adaptations. As discussed in Section \ref{sec:noise}, this noise can be used to emulate the wander of a bad quality clock.

We used different values of $g_{ij}^w$ and $g_i^d$ to differentiate the noise conditions in the optimization scheme. The large jitter scenario is represented in by $g_{i}^d=1e-3$ $\forall i$, $g_{21}^w=g_{31}^w=100$ and $g_{ij}^w=1$ otherwise. The large wander scenario is represented by $g_{i}^d=1e-1$ $\forall i$ and $g_{ij}^w=1$. Finally, the large jitter and wander scenario is represented using $g_{i}^d=1e-1$ $\forall i$, $g_{21}^w=g_{31}^w=100$ and $g_{ij}^w=1$ otherwise. The output parameter values for all three cases are present also in Figure \ref{fig:exp_h2}.

\begin{figure}[htp]
\centering
\includegraphics[width=1.0\columnwidth]{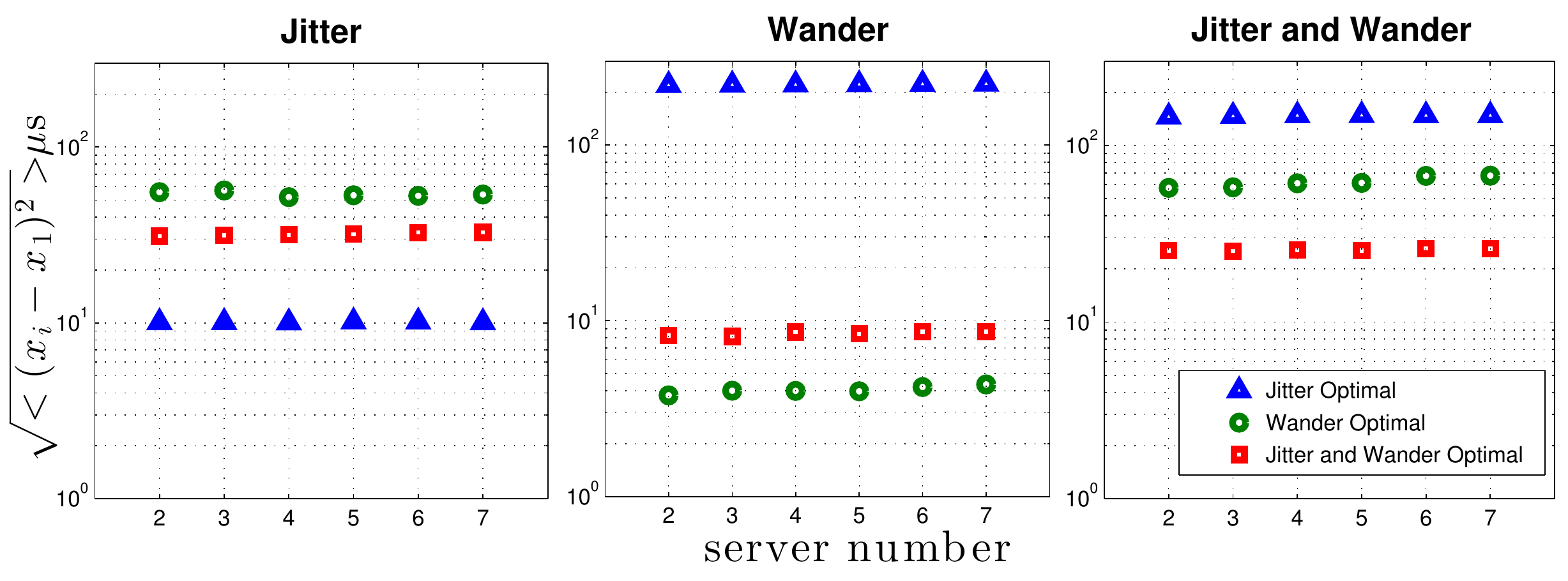}
\caption{$\mathcal H_2$ Performance optimization: offset variance vs server number}\label{fig:H2_cases}
\end{figure}

Figure \ref{fig:H2_cases} shows the  standard deviation of the offset between servers serv2-serv7 and serv1 in the three experimental scenarios and for the three different set of parameters shown in Figure \ref{fig:exp_h2}. It can be seen that although the configuration tuned for jitter performs very well in cases with large jitter, it performs quite poorly in scenarios with large wander. Similarly, the configuration tuned for wander does not perform well in high jitter scenarios.
However, the configuration tuned for jitter and wander is able to provide acceptable performance in all three experimental scenarios. Thus, we experimentally demonstrate a fundamental tradeoff between jitter and wander.

\changed{Finally, it is interesting to notice that due to the fact the optimization is solved using $v_i(t_k)=x_i(t_k)-x_1(t_k)$ as performance metric, 
the choice of parameters does not degrades the performance of each clock with the hop count.}

\section{Conclusion}\label{sec:conclusions}
This paper presents a clock synchronization protocol that is able to synchronize networked nodes without explicit estimation of the clock skews and steep corrections on the time. Our solution is guaranteed to converge even in the presence of timing loops which allow different clients to share timing information and even collectively outperform individual clients when the time source has large jitter.
The system is robust to noisy measurements and wander provided that the topology has a well defined leader, and we can optimize the parameter values to minimize noise variance. \changed{We implemented our solution on a cluster of IBM BladeCenter servers and empirically verified our predictions and our protocol's supremacy over some existing solutions.}

\changed{Further evaluation of our protocol is needed. In particular, we are interested in comparing our solution with other protocols such as PTP and RADclock, as well as studying its robustness under stressed scenarios. Another interesting direction is to devise a distributed algorithm, exploiting our optimization framework, that can adapt the parameter values depending on the network condition.}


%
%

\bibliography{biblio}
\bibliographystyle{IEEEtran}

\appendices

\section{ Proof of Lemmas }\label{app:lemmas}

\subsection{ Proof of Lemma \ref{lem:multiplicity} }\label{app:lemma1}
\begin{IEEEproof}
We first compute the characteristic polynomial
\begin{align*}
&\det(\lambda I_{3n}-A) =
\left|
\begin{array}{ccc}
(\lambda-1)I_n  & -\tau R &\0_{n\times n} \\
\kappa_1L & (\lambda-1)I_n & \kappa_2I_n \\
pL & 0 & (\lambda-1+p) I_n
\end{array}
\right| \\
&=(\lambda-1)^n
\left|
\begin{array}{cc}
(\lambda-1)I_n+\frac{\tau \kappa_1}{\lambda-1}LR &  \kappa_2I_n\\
\frac{\tau p}{\lambda-1} LR & (\lambda -1 +p) I_n
\end{array}
\right|\\
&=\det\left( (\lambda-1)^2(\lambda-1+p)I_n +  [(\lambda-1)\kappa_1 \right.  \\
&\left.  +\changed{p(\kappa_1-\kappa_2)}]\tau LR\right) = \prod_{l=1}^{n} g_l(\lambda),
\end{align*}
where $g_l(\lambda)$ is as defined in \eqref{eq:g_l} and we have  iteratively use the determinant property of block matrices
$
\det(A) = \det(A_{11}) \det(A\backslash A_{11})
$
where
$
A=\left[\begin{array}{cc}
A_{11} & A_{12}\\
A_{21} & A_{22}
\end{array}
\right]
$ and $A\backslash A_{11}=A_{22} - A_{21}A_{11}^{-1}A_{12}$ is the Schur complement of $A_{11}$~\cite{horn_matrix_1990}.


Thus, $\lambda=1$ is a double root of the characteristic polynomial if and only if $\kp_1\neq\kp_2$, $p>0$ and $\tau LR$ has a simple zero eigenvalue, i.e. \eqref{eq:nu_condition}.
Now, since $R$ is nonsingular \eqref{eq:nu_condition} must hold for the eigenvalues of $L$ as well, which is in fact true if and only if the directed graph $G(V,E)$ is connected~\cite{xie_consensus_2011}.
\end{IEEEproof}

\subsection{Proof of Lemma \ref{lem:jordan_chains} }\label{app:lemma2}

\begin{IEEEproof}
We start by computing the right Jordan chain. By definition of $\zeta_1$,
$
(A - I)\zeta_1 = 0_n.
$
Thus, if $\zeta_1 = [x_1^T\;s_1^T\;y_1^T]^T$, then the following system of equations must be satisfied
\begin{align}
\tau &R s_1 = \0_n\text{ (a), } -\kappa_1 Lx_1  - \kappa_2 y_1 =\0_n\text{ (b) }\text{ and }\nonumber\\
 -p &Lx_1 - p y_1 =\0_n\text{ (c).} \label{eq:eq_sys}
\end{align}
Equation (\ref{eq:eq_sys}a) implies $s_1=0$. Now, since $p>0$, (\ref{eq:eq_sys}c) implies $Lx_1=-y_1$, which when substituted in (\ref{eq:eq_sys}b) gives
$
(\kappa_2-\kappa_1)y_1=\0_n.
$
Thus, since $\kappa_1\neq\kappa_2$, $y_1=\0_n$ and $x_1\in\ker(L)$. By choosing $x_1=\alpha_1\1_n$ (for some $\alpha_1\neq 0$) we obtain
$
\zeta_1= \alpha_1
\left[
\1_n^T \; \0_n^T \; \0_n^T
\right]^T.
$

Notice that the computation also shows that $\zeta_1$ is the unique eigenvector of $\mu(A)=1$ which implies that there is only one Jordan block, of size 2. The second member of the chain, $\zeta_2$, and $\zeta_3$ can be computed similarly by solving
$(A - I_n)\zeta_2 = \zeta_1$ and $(A-(1-p)I_n)\zeta_3=\0_n$. This gives
\begin{align*}
\zeta_2= \left[\begin{array}{c}
\alpha_2\1_n \\ \frac{\alpha_1}{\tau}R^{-1}\1_n \\ \0_n
\end{array}\right]
\quad\text{ and }\quad
\zeta_3= \alpha_3\left[\begin{array}{c}
-\frac{\tau\kp_2}{p^2}\1_n \\ \frac{\kp_2}{p}R^{-1}\1_n \\ R^{-1}\1_n
\end{array}\right].
\end{align*}
In computing $\zeta_3= [x_3^T\;s_3^T\;y_3^T]^T$, we obtain $Lx_3=0$ and $Rx_3=-\frac{\tau}{p}s_3=-\frac{\kp_2\tau}{p^2}y_3$. $\zeta_3$ follows by taking $y_3=\alpha_3R^{-1}\1_n$.

The vectors $\eta_1$, $\eta_2$ and $\eta_3$ can be solved in the same way using $\eta_2^T(A - I)=\0_n^T$,
$\eta_1^T(A - I)=\eta_2^T$ and $\eta_3^T(A - (1-p)I)=\0_n^T$. This gives
$\eta_1= \left[
\frac{\beta_2}{\tau}R^{-1}\xi^T \; \beta_1\xi^T \; (-\frac{\kp_2}{p}\beta_1 + \frac{\kp_2}{p^2}\beta_2)\xi^T
\right]^T
$,
$\eta_2= \beta_2\left[
\0_n^T \; \xi^T \; \frac{\kp_2}{p}\xi^T
\right]^T
\quad\text{ and }\quad
\eta_3= \beta_3\left[
\0_n^T \; \0_n^T \; \xi^T
\right]^T.
$
We set  $\alpha_1=\alpha_2=\alpha_3=1$; this can be done without loss of generality provided we still satisfy $\eta_l^T\zeta_l=1$ and $\eta_l^T\zeta_h=0$ for $l\neq h$.
\changed{Finally, $\eta_1^T\zeta_1=1$ gives $\beta_2=\gamma\tau$, $\eta_3^T\zeta_3=1$ gives $\beta_3=\gamma$ and $\eta_1^T\zeta_2=0$ gives $\beta_1=-\beta_2=-\gamma\tau$.}
\end{IEEEproof}

%
%

\section{Proof of Theorem \ref{th:convergence} }\label{app:th:convergence}
\begin{IEEEproof}\\
\changed{
\noindent{\it 1) $\implies$ 2):} Since we are under the conditions of Lemmas 1 and 2, then we can use \eqref{eq:hatA} and since $\rho(\hat J_2)<1$  all the eigenvalues of $\hat A$ are within the unit circle, i.e. $\rho(\hat A)<1$. Therefore, it follows that $\delta z_k =\hat A^k z_0 \rightarrow \0_{3n\times 3n}z_0=\0_{3n}$. To show \eqref{eq:tildez_sync} we first notice that 
\begin{subequations}\label{eq:tilde}
\begin{align}
\tilde x(t_{k+1})&= \tilde x(t_k) + \tau \tilde s(t_k)\label{eq:tildex}\\
\tilde s(t_{k+1})&= \tilde s(t_k) - \kappa \tilde y(t_k)\label{eq:tildes}\\
\tilde y(t_{k+1})&= (1-p) \tilde y(t_k)\label{eq:tildey}
\end{align}
\end{subequations}
Therefore, since $|1-p|<1$, \eqref{eq:tildey} implies that $\tilde y(t_k)\rightarrow 0$. Thus, by \eqref{eq:tildes} we also have $\tilde s(t_k)\rightarrow s^*$ for some $s^*$, which also implies that $\tilde x(t_{k+1})- \tilde x(t_k)\rightarrow \tau s^*$ giving $\tilde x(t_k)\rightarrow x^\text{ref}(t_k)=x^* + (t_k-t_0)s^*$ for some $x^*$.
}

\noindent\changed{
\noindent{\it 2) $\implies$ 3):} This follows directly from \eqref{eq:average} and \eqref{eq:deviations}.
}

\noindent\changed{
\noindent{\it 3) $\implies$ 1):} The algorithm achieves synchronization whenever \eqref{eq:synchronization} holds. Then, it follows from \eqref{eq:system_z} and \eqref{eq:synchronization} that asymptotically the system behaves according to
\begin{align*}
z_{k} &= \left[
\begin{array}{c}
x_{k} \\ s_{k}\\ y_{k}
\end{array}
\right] =
\left[
\begin{array}{c}
x^*\1_n \\ r^* R^{-1}\1_n\\ \0_n
\end{array}
\right] + k\left[
\begin{array}{c}
\tau r^* \1_n \\ \0_n\\ \0_n
\end{array}
\right] \\
&=  \left(\tau r^*\zeta_2 + (x^*-\tau r^*)\zeta_1\right) + kr^*\tau \zeta_2.
\end{align*}
\changed{Thus, since $P:=[ \zeta_1 \quad ... \quad\zeta_{3n}]$ is invertible, its columns $\zeta_l$ are linearly independent.} Therefore, if the system synchronizes for arbitrary initial condition, then it must be the case that the effect of the remaining modes $\mu_l(A)$ vanishes, which can only happen if for every $\mu_l(A)\neq 1$, $|\mu_l(A)|<1$ and the multiplicity of $\mu_l(A)=1$ is two, i.e. $\rho(\hat J_2)<1$.. Now suppose that either $G(V,E)$ is not connected, $\kappa_1=\kappa_2$, $p=0$. Then by Lemma \ref{lem:multiplicity}, the multiplicity of $\mu_l(A)=1$ is not two which is a contradiction. Similarly, if $p>2$, $p<0$ then the system has at least one eigenvalue $|\mu_l(A)|>1$. 
 Thus, we must have $\rho(\hat J_2)<1$, $\kappa_1\neq \kappa_2$, $2>p>0$ and $G$ connected whenever the system synchronizes for arbitrary initial condition.}

Finally, to obtain \eqref{eq:xs&rs} we use a similar computation to the one of Lemma 1 to show that $\tilde P$ and $\tilde P^{-1}$ in \eqref{eq:tildeA} are given by
\begin{align*}
\tilde P  
= 
\left[
\begin{array}{ccc} 
 1 & 1 & -\tau\frac{\kappa_2}{p^2} \\
  0& \frac{1}{\tau} & \frac{\kp_2}{p} \\
  0 & 0 & 1
\end{array}
\right],
\tilde P^{-1} 
= 
\left[
\begin{array}{ccc} 
1 & -\tau  & (\frac{1}{p^2}+\frac{1}{p})\kp_2\tau \\
0  & \tau & -\tau\frac{\kp_2}{p}\\
 0 & 0 & 1
\end{array}
\right].
\end{align*}
Thus, since $\tilde A^k = \tilde P \hat J_1^k \tilde P^{-1}$ a direct computation shows that
\[
\tilde A^k = \tilde \zeta_1\tilde\eta_1^T + \tilde \zeta_2\tilde\eta_2^T + k\tilde \zeta_1\tilde\eta_2^T+(1-p)^k\tilde \zeta_3\tilde\eta_3^T
\]
where $\tilde P = [\tilde \zeta_1 \tilde \zeta_2 \tilde \zeta_3]$ and $\tilde P^{-1} = [\tilde \eta_1 \tilde \eta_2 \tilde \eta_3]^T$
Therefore,
\[
\tilde z_k = \tilde A^k\tilde z_0\rightarrow
\left[
\begin{array}{ccc}
1 & \tau k & -\tau k \frac{\kp_2}{p} + \tau\frac{\kp_2}{p^2}\\
0 & 1 & -\frac{\kp_2}{p}\\
0 & 0 & 0
\end{array}
\right]\left[\begin{array}{c} \tilde x(t_0)\\\tilde s(t_0)\\\tilde y(t_0)\end{array}\right]
\]
which implies that $\tilde y(t_k)\rightarrow 0$, $\tilde s(t_k)\rightarrow  \tilde s(t_0) -\frac{\kp_2}{p}y(t_0)$ and 
$\tilde x(t_k) \rightarrow \tilde x(t_0) +  \tau\frac{\kp_2}{p^2}\tilde y(t_0)+ \tau k (s(t_0) -\frac{\kp_2}{p}y(t_0))$. Result follows by definition of $x^*$ and $r^*$ in \eqref{eq:xs&rs} and definition of $x^\text{ref}(t_k)$ in \eqref{eq:synchronization}.

\end{IEEEproof}

\section{ Proof of Theorem \ref{th:param_sync} }\label{app:th:param_sync}

\begin{IEEEproof}
We will show that when $G(V,E)$ is connected with $\mu(L)\in\mathds R$, then (i)-(iii) are equivalent to the conditions of Theorem \ref{th:convergence}.

Since, $G(V,E)$ is connected and (i)-(ii) satisfies $p>0$ and $\kp_1\neq\kp_2$, the conditions of  Lemma \ref{lem:multiplicity} are satisfied. Therefore the multiplicity of $\mu(A)=1$ is two and by \eqref{eq:nu_condition} these are the roots of
$g_n(\lambda) = (\lambda-1)^2(\lambda -1 +p),$
 which corresponds to the case $\nu_n=0$.
Thus, to satisfy Theorem \ref{th:convergence} we need to show that the remaining eigenvalues are strictly in the unit circle. This is true for the remaining root of $g_n(\lambda)$ if and only if (i).


For the remaining $g_l(\lambda)$, this implies that are Schur polynomials. Thus, we will show that $g_l(\lambda)$ is a Schur polynomial if and only if (i)-(iii) hold.
We drop the subindex $l$ for the rest of the proof.

We first transform the Schur stability problem into a Hurwitz stability problem. Consider the change of variable $\lambda = \frac{s+1}{s-1}$. Then $|\lambda|<1$ if and only if $\mathds R[s]<0$.

%

Now, since $\nu>0$ by \eqref{eq:nu_condition}, let
\begin{align*}
P(s) & =\frac{(s-1)^3 }{\delta\kp p\nu}g\left(\frac{s+1}{s-1}\right) =
 s^3 + \left(\frac{2\kappa_1}{\delta\kp p}-3\right)s^2 \\&+ \left( \frac{4}{\delta\kp\nu} + 3 -\frac{4\kappa_1}{\delta\kp p} \right)s + \frac{4(2-p)}{\delta\kp p\nu} +\frac{2\kp_1}{\delta\kp p}-1
\end{align*}
where $\delta\kp = \kp_1 -\kp_2$.

We will apply Hermite-Biehler Theorem  to $P(s)$, but first let us express what $1)$ and $2)$ of the Theorem mean here.

Condition $1)$  becomes
\begin{align}\label{eq:HB1}
\frac{2\kappa_1}{\delta\kp p}-3 > 0.
\end{align}

Now let $P^r(\omega)$ and $P^i(\omega)$ be as in Hermite-Biehler Theorem, i.e. let
\begin{align*}
P^r(\omega)=& -\left(\frac{2\kappa_1}{\delta\kp p}-3\right)\omega^2+\frac{4(2-p)}{\delta\kp p\nu} +\frac{2\kp_1}{\delta\kp p}-1\\
P^i(\omega)=&-\omega^3 + \left( \frac{4}{\delta\kp\nu} + 3 -\frac{4\kappa_1}{\delta\kp p} \right)\omega
\end{align*}

The roots of $P^r(\omega)$ and $P^i(\omega)$ are given by $\omega_0=\pm\sqrt{\omega^r_0}$ and $\omega_0\in\{0,\;\pm\sqrt{\omega^i_0}\}$ respectively, where
\begin{equation}\label{eq:root}
\omega^r_0:=\frac{\frac{4(2-p)}{\delta\kp p\nu} +\frac{2\kp_1}{\delta\kp p}-1}{\frac{2\kappa_1}{\delta\kp p}-3}\text{ and }
\omega^i_0:= \frac{4}{\delta\kp\nu} + 3 -\frac{4\kappa_1}{\delta\kp p}
\end{equation}

Since the roots $P^r(\omega)$ and $P^i(\omega)$ must be real, we must have $\omega_0^r>0$ and $\omega_0^i>0$. Therefore, by monotonicity of the square root, the interlacing condition $2)$ is equivalent to
\begin{equation}\label{eq:HB2}
0<\omega_0^r<\omega_0^i.
\end{equation}
Thus we will show:  (i)-(iii) hold $\iff$ \eqref{eq:HB1} and \eqref{eq:HB2} hold.

It is straightforward to see that using (i) and (ii) we can get \eqref{eq:HB1}. On the other hand,  $\omega_o^i>0$ from \eqref{eq:HB2} together with \eqref{eq:HB1} gives
$
0<\frac{4}{\delta\kp\nu} + 3 -\frac{4\kappa_1}{\delta\kp p} <\frac{4}{\delta\kp\nu}
$, which implies that $\delta\kp>0$, and therefore (ii) follows.

Now using \eqref{eq:HB1} and the definition of $\omega^r_0$ in \eqref{eq:root}, $\omega_0^r>0$ becomes
$
\frac{4(2-p)}{\delta\kp p\nu} +\frac{2\kp_1}{\delta\kp p}-1>0
$
which always holds under (i) and (ii) since the first term is always positive and
$
\frac{2\kp_1}{\delta\kp p}-1>\frac{2\kp_1}{\delta\kp p}-3>0
$
by \eqref{eq:HB1}.

Using \eqref{eq:root}, $\omega_0^r<\omega_0^i$ is equivalent to
%
\begin{equation}
\nu < \frac{p(\kp_2 - \delta\kp p)}{(\kp_1-\delta\kp p)^2}. \label{eq:inequality}
\end{equation}
Finally, $\nu_l=\mu_l(\tau LR)=\tau \mu_{l}(LR)$. Thus, since \eqref{eq:inequality} should hold $\forall l\in\{1,...,n-1\}$, then
\[
\tau < \min_l \frac{p(\kp_2 - \delta\kp p)}{\mu_l(LR)(\kp_1-\delta\kp p)^2} = \frac{p(\kp_2 - \delta\kp p)}{\mu_{\max}(\kp_1-\delta\kp p)^2}
\]
which is exactly (iii).
\end{IEEEproof}

\section{Proof of Theorems \ref{th:frequency_drift} and  \ref{th:time_offsets}} \label{app:th:frequency_drift_and_time_offsets}

\noindent{\it Proof of Theorem \ref{th:frequency_drift}:}\\
Using \eqref{eq:deviations}, \eqref{eq:system_noise} and \eqref{eq:Bs} we can modify \eqref{eq:tilde} to get
\begin{subequations}\label{eq:btilde}
\begin{align}
\tilde x(t_{k+1})&= \tilde x(t_k) + \tau \tilde s(t_k)\label{eq:btildex}\\
\tilde s(t_{k+1})&= \tilde s(t_k) - \kappa_2 \tilde y(t_k) + \kp_1 \tilde w \label{eq:btildes}\\
\tilde y(t_{k+1})&= (1-p) \tilde y(t_k) + p \tilde w\label{eq:btildey}
\end{align}
\end{subequations}
where 
$
\tilde w = -\xi^TB_G^{-}\diag[\alpha_{ij}g_{ij}^w]\bar w = \sum_{i=1}^{n} \xi_i \sum_{j\in\mathcal N_i} \alpha_{ij}g_{ij}^w\bar w_{ij}.
$ 
It follows then that \eqref{eq:btildey} implies that $\tilde y(t_k)\rightarrow \tilde w$ which implies that 
$
\tilde s(t_{k+1})- \tilde s(t_{k})\rightarrow (\kp_1 - \kp_2)\tilde w.
$

Therefore, since $\kp_1\neq\kp_2$, $\tilde s(t_k)$ constantly drifts unless 
\begin{equation}\label{eq:null_freq_drift_condition}
\tilde w =-\xi^TB_G^{-}\diag[\alpha_{ij}g_{ij}^w]\bar w =0.
\end{equation}
Finally,  there are two different scenarios in which \eqref{eq:null_freq_drift_condition} can be satisfied.
\begin{enumerate}
\item $G$ has a  unique leader (say $i=1$): In this case we have $\mathcal N_1=\emptyset$, i.e. $\alpha_{1j}=0$ $\forall j$,  $\xi_1=1$ and $\xi_j=0$ $\forall j\neq 1$. That is $-\xi^TB_G^{-}\diag[\alpha_{ij}g_{ij}^w]\bar w=\xi_10=0$
\item $G$  does not have a well defined root: Thus, there are at least two nodes with $\xi_i\neq 0$ and $\bar w$ is such that $\xi^TB_G^{-}\diag[\alpha_{ij}g_{ij}^w]\bar w=0$.
\end{enumerate}
However, $2)$ is only satisfied by a set of values of $\bar w$ with zero measure. Thus, there should be a unique leader for synchronization. \qed
\vspace{.2cm}

\noindent{\it Proof of Theorem \ref{th:time_offsets}:}\\
Similar to the proof of Theorem \ref{th:frequency_drift}, the evolution of $\delta z_k$ can be described using
\begin{subequations}\label{eq:delz}
\begin{align}
&\delta x_{k+1} = \delta x_k + \tau R\delta s_k\label{eq:delx}\\
&\delta s_{k+1} = -\kp_1 L \delta x_k + \delta s_k -\kp_2 \delta y_k +\kp_1\delta w\label{eq:dels}\\
&\delta y_{k+1} = -p L\delta x_k + (1-p)\delta y_k +p\delta w\label{eq:dely}
\end{align}
\end{subequations}
%

Now, since $\rho(NA)<1$, then $\delta z_k\rightarrow \delta z^*$, where $\delta z^*$ is a fixed point of \eqref{eq:delz}.
Thus, \eqref{eq:delx} implies that $\delta\bar s^*=0$ and \eqref{eq:dels}$-\frac{\kp_1}{p}$\eqref{eq:dely} gives
$
(\kp_1-\kp_2)\delta\bar y^* = 0,
$
which implies $\delta\bar y^*=0$ since $\kp_1>\kp_2$.
Finally using \eqref{eq:dely} again we have
\begin{align}
L\delta\bar x^* &= \delta w \\
L^\dagger L \delta\bar x^* &= L^\dagger  \delta w \label{eq:step2}\\
N_1N_3 \delta\bar x^*& = N_1L^\dagger \delta w\label{eq:step4}\\
\delta\bar x^* &= N_1L^\dagger \delta w\label{eq:step5}
\end{align}
where in \eqref{eq:step2} we multiplied by $L^\dagger$, in \eqref{eq:step4} we used $N_3 := L^\dagger L=(I_n-\frac{1}{n}\1_n\1_n^T)$ and left multiply by $N_1$, and in \eqref{eq:step5} we used de identities $N_1N_3=N_1$ and $N_1\delta\bar x=N_1^2 \bar x=N_1\bar x=\delta\bar x$.
\qed


\section{$\mathcal H_2$ Optimization using Hifood} \label{app:Hifood}

The software package Hifood~\cite{popov_fixed-structure_2010} does not solve \eqref{eq:H2} directly. Instead, it solves:
\begin{subequations}\label{eq:H22}
\begin{align}
&\quad\min_{K,X} \;\;\; f(K):=\sqrt{\trace [X\bar B\bar B^T]} \\
&\quad\text{ subject to } \quad \rho(\bar A  )\leq \rho^* \\
& \quad X=\bar A^TX\bar A + \bar C^T\bar C
\end{align}
\end{subequations}
where $\bar A := A_1 +  B_2K C_2$, $\bar B:= B_1 + B_2 K  D_{21}$ and $\bar C:= C_1$. In this formulation $\delta z_k$ is interpreted as evolving according to the closed loop standard form system
\begin{align*}
 \delta z_{k+1} &=  ( A_1 +  B_2K C_2)\delta z_k + ( B_1 + B_2 K  D_{21}) e_k \\
 v_{k} &=  C_1 \delta z_k ,
\end{align*}
and the optimization variable $K$ is the static-output feedback matrix.

Therefore, to use Hifood we first need to rewrite \eqref{eq:H2} using \eqref{eq:H22}.
This can be done by setting
%
%
{\small
\begin{align*}
& A_1 = \hat A, \;  C_1= \hat C,  C_2 = \left[\begin{array}{ccc} B_G^T & \0_{m\times n}& \0_{m\times n}\\ \0_{n\times n} & I_n & \0_{n\times n}\\ \0_{n\times n}& \0_{n\times n} & I_n \end{array}\right],\\
& B_2 = \left[\begin{array}{ccccc} N_1 R& \0_{n\times m} & \0_{n\times n} &\0_{n\times m} & \0_{n\times n}\\
\0_{n\times n}& B_G^{-} & N_2 & \0_{n\times m} & \0_{n\times n}\\
							\0_{n\times n}& \0_{n\times m} & \0_{n\times n} &B_G^{-} & N_2
\end{array}
\right],\\
& B_1 = \left[ \begin{array}{cc}
\0_{n\times m} & \0_{n\times n}\\
\0_{n\times m} & \diag[g_i^d]\\
\0_{n\times m} & \0_{n\times n}
\end{array}\right],\;
  D_{21} = \left[\begin{array}{ccc} \diag[g_{ij}^w] & \0_{m\times n} \\ \0_{n\times m} & \0_{n\times n} \\ \0_{n\times m} & \0_{n\times n}
\end{array}\right],\allowbreak\\
&\text{ and } K= \left[\begin{array}{ccc} \0_{n\times m} & \tau I_n & \0_{n\times n}\\ -\kp_1 \diag[\alpha_{ij}] & \0_{m\times n} & \0_{m\times n} \\ \0_{n\times m} & \0_{n\times n} & -\kp_2 I_n \\
 -p \diag[\alpha_{ij}] & \0_{m\times n} & \0_{m\times n} \\ \0_{n\times m} & \0_{n\times n} & -p I_n \\
\end{array}\right].
\end{align*}
}
Using these definitions it is straight forward to verify that $( A_1+  B_2 K  C_2) = \hat A$, $ B_1 +  B_2 K  D_{21}=\hat B$ and $ C_1 = \hat C$.

The main difficulty in solving \eqref{eq:H2} instead of \eqref{eq:H22} is that our controller $K$ is a nonlinear function of the parameters $K(\kp_1,\kp_2,p,\alpha)$ and cannot be readily obtained using \eqref{eq:H22}.
Furthermore, the main source of nonlinearity comes from the products $\kp_1\diag[\alpha_{ij}]$ and $p\diag[\alpha_{ij}]$. This structure is not currently supported by traditional software distributions, which usually only support sparsity patterns, and therefore needs to be implemented.

Fortunately, Hifood only uses gradient information in their implementation of BGS and gradient bundle stages. Thus, to implement discrete time $\mathcal H_2$ optimization we generated a new Matlab subroutine that evaluated the $\mathcal H_2$ norm $f$ as well as  its gradients.

The evaluation of the gradient is performed in three stages using the chain rule. We first compute the gradients of $f$ with respect to 
$\bar A$, $\bar B$ and $\bar C$ which are given by
\[
\nabla_{\bar A}f=\frac{1}{f}X\bar A Y, \quad \nabla_{\bar B}f=\frac{1}{f}X\bar B\;\;\;\text{   and   }\;\;\;\nabla_{\bar C}f=\frac{1}{f}\bar C Y
\]
\changed{where $Y$ is the solution to $Y=\bar AY\bar A^T + \bar B\bar B^T$.}

Once $\nabla_{\bar A} f$, $\nabla_{\bar B} f$ and $\nabla_{\bar C} f$ are computed we can use the subroutines of hifood to compute $\frac{\partial \bar A}{\partial K}$, $\frac{\partial \bar B}{\partial K}$ and $\frac{\partial \bar C}{\partial K}$. Finally, we obtain
{\small\begin{align*}
\nabla_{\kp_1}f=\trace \left[\left(\nabla_{\bar A}f^T\frac{\partial \bar A}{\partial K} + \nabla_{\bar B}f^T\frac{\partial \bar B}{\partial K} +
			\nabla_{\bar C}f^T\frac{\partial \bar C}{\partial K}  \right)\frac{\partial K}{\partial \kp_1}\right]
\end{align*}}
and similarly for other parameters.

\end{document}